\documentclass[leqno]{mnarticle} %

\usepackage{amssymb}
\usepackage{graphicx} 
\usepackage{xspace}
\usepackage{multirow}
\usepackage{bbm}
\usepackage{Macros/mydef} 
\usepackage{algpseudocode}
\usepackage[numbers]{natbib}
\let\cite=\citet
\usepackage{enumerate}
\usepackage{booktabs}
\usepackage{graphicx}
\usepackage{subcaption}
\usepackage{todonotes}
\usepackage{thmtools}
\usepackage{cleveref}

\usepackage{graphicx}
\usepackage{subcaption}
\usepackage{array}
\usepackage{graphicx}
\usepackage{booktabs}

\ifpdf
  \DeclareGraphicsExtensions{.eps,.pdf,.png,.jpg}
\else
  \DeclareGraphicsExtensions{.eps}
\fi

\newsiamremark{remark}{Remark}
\newsiamremark{hypothesis}{Hypothesis}
\crefname{hypothesis}{Hypothesis}{Hypotheses}
\newsiamthm{claim}{Claim}
\newsiamthm{assumption}{Assumption}

\usepackage{amsopn}

\providecommand{\bsfa}{\textsf{\textbf a}}
\providecommand{\bsfr}{\textsf{\textbf r}}
\providecommand{\bsfs}{\textsf{\textbf s}}
\providecommand{\bsfq}{\textsf{\textbf q}}
\providecommand{\bsfp}{\textsf{\textbf p}}
\providecommand{\bsfb}{\textsf{\textbf b}}
\providecommand{\bsfd}{\textsf{\textbf d}}
\providecommand{\bsfj}{\textsf{\textbf j}}

\providecommand{\bsfc}{\textsf{\textbf c}}

\providecommand{\bsft}{\textsf{\textbf t}}
\providecommand{\bsfP}{\textsf{\textbf P}}
\providecommand{\bsfV}{\textsf{\textbf V}}
\providecommand{\Mone}{M_1}
\providecommand{\Mk}{M_k}
\newcommand{\Vf}{\calV_{\rm f}}

\newcommand{\Vd}{\calV_{\rm D}}

\newcommand{\Ef}{\calE_{\rm f}}
\newcommand{\Nf}{\calN_{\rm f}}
\newcommand{\mk}{m^{\polP_k}}
\newcommand{\mone}{m^{\polP_1}}
\newcommand{\IoneF}{I_h^{\polP_1}\polf^n}
\newcommand{\Fh}{\polf_h^n}
\newcommand{\dH}{\bsfd^{{\rm H},n+1}}
\newcommand{\Un}{\bsfU^{n}}
\newcommand{\UH}{\bsfU^{{\rm H},n+1}}
\newcommand{\UL}{\bsfU^{{\rm L},n+1}}
\newcommand{\ULhat}{\widehat{\ \bsfU \ }^{{\rm L}, n+1}}
\newcommand{\ULtilde}{\widetilde{\ \bsfU \ }^{{\rm L}, n+1}}
\newcommand{\Ubar}{\overline{\bsfU}^n}
\newcommand{\V}{\bsfV^{n+1}}

\usepackage{mathabx}

\usepackage{etoc}

\headers{Graph-Poisson limiting}{Murtazo Nazarov}

\title{High-order invariant-domain preserving continuous finite elements for the compressible Euler equations via graph-Poisson flux reconstruction
  \thanks{\funding{This research is funded by Swedish Research Council (VR) under Grant Number 2021-04620.}}}

\author{
Murtazo Nazarov\thanks{
Division of Scientific Computing, Department of Information Technology,
Uppsala University, Uppsala 751 05, Sweden
(\email{murtazo.nazarov@uu.se}).
}
}

\providecommand{\tdel}[1]{}
\providecommand{\mdel}[1]{}
\providecommand{\draftnote}[1]{}

\begin{document}

\maketitle

\begin{abstract}
  We propose a high-order invariant-domain preserving (IDP) continuous finite
  element method for the compressible Euler equations, with Lagrange elements of
  degree one to four on unstructured triangular meshes. A stabilized high-order
  Galerkin method is combined with a low-order IDP graph-viscosity method on the
  fine $\polP_1$ mesh whose vertices are the nodes of the high-order mesh. The
  difference between the two updates, each with its own lumped mass, is written
  as antisymmetric edge fluxes by a graph-Poisson flux reconstruction, either
  globally or element by element; the element-local version needs no global
  linear solve. Convex limiting of these fluxes gives a conservative method that
  keeps local bounds on the density, the positivity of the internal energy and a
  local minimum principle on the specific entropy, also for elements with zero or
  negative lumped masses such as $\polP_2$ and $\polP_4$; for these elements,
  the high-order mass is restored within the bounds by an element-local exchange,
  completed if needed by a global round. Numerical benchmarks
  show near-optimal convergence rates for smooth solutions and robustness for
  problems with strong shocks.
\end{abstract}



\section{Introduction}
\label{Sec:introduction}
In this paper, we approximate the system of compressible Euler equations
with a continuous finite element method. Our aim is a method that, for Lagrange elements of any polynomial degree, preserves the so-called {\em invariant domain property}, see \eg
\citep{Chueh_Conley_Smoller_1977} and
\citep[Sec.~2]{Guermond_Popov_2016}. That is, we seek a
high-order invariant-domain preserving (IDP) finite element approximation of the Euler equations
in the fully discrete setting.

Most high-order IDP methods consist of three ingredients:
an IDP low-order scheme, a high-order candidate, and a limiter that
blends them under quasiconcave constraints. This approach is
usually referred to as the {\em flux-corrected transport} (FCT) paradigm,
see for instance \citep{Boris_Book_1973, Zalesak_1979, Kuzmin_Lohner_Turek_2012,
Kuzmin_2002, Guermond_Nazarov_Popov_Yong_2014} and references therein. In FCT, the high-order update is written as the low-order update plus
antidiffusive edge fluxes, which are cut back until the bounds hold.
For the compressible Euler
equations, the constraints are positivity of the density and of the internal energy, and a minimum principle
on the specific entropy. First-order IDP continuous finite elements
for hyperbolic systems, based on a {\em graph viscosity} with the maximum
wave speed of local Riemann problems, are presented in \citep{Guermond_Popov_2016,
  Guermond_Popov_Fast_Riemann_2016}.
A higher-order extension
of the first-order method, with the entropy viscosity stabilization of
\citep{Guermond_pasquetti_popov_JCP_2011}
as the high-order candidate and a {\em convex limiting} algorithm, was later
proposed in \citep{Guermond_Nazarov_Popov_Tomas_2018}.
More generally, the high-order candidate can be any conservative method whose
difference from the low-order update can be written as antisymmetric fluxes on
the edges of the low-order stencil, for instance with the entropy or residual
viscosity stabilization of \citep{Guermond_pasquetti_popov_JCP_2011, Nazarov_2013}.
A
monolithic version for continuous finite elements, in which the limited fluxes
enter the semi-discrete scheme, was introduced in \citep{Kuzmin_2020}.

The method of \citep{Guermond_Nazarov_Popov_Tomas_2018} was shown to be robust
for various challenging shock problems; however, it
is restricted to linear elements and is at most second-order
accurate. There are two main obstacles to extending this approach
to higher polynomial degrees. First, as shown in \citep{Guermond_Popov_Yang_2017},
the consistent mass matrix breaks the invariant domain property; therefore,
the low-order IDP scheme needs a lumped mass with positive entries.
Unfortunately, this is not the case for many Lagrange elements on triangles:
for example, quadratic polynomials $\polP_2$ produce zero lumped mass entries at the vertices,
while quartic polynomials $\polP_4$ give some negative lumped masses.
The second obstacle is the amount of low-order viscosity on high-order stencils:
the time-step bound drops fast and the diffusion grows with the polynomial degree,
see \citep[Sec.~4.1]{Guermond_Nazarov_Popov_2024} and references therein.

Getting around these obstacles has been one of the main directions in
the development of IDP schemes for conservation laws in recent years.
One option is to use Bernstein elements, whose lumped masses are positive, see \citep{Lohmann_etal_2017, Hajduk_etal_2020a, Hajduk_etal_2020b}; in
\citep{Kuzmin_Hajduk_Vedral_2026}, a Bernstein discretization with a nonlinear
stabilization is limited element by element towards intermediate cell averages.
Another approach, inspired by the high-order continuous residual distribution schemes of
\citep{Abgrall_etal_2017} and \citep{Abgrall_Bacigaluppi_Tokareva_2019}, is to modify the low-order fluxes so that they carry the high-order
mass. This idea was later
extended to continuous finite element discretizations with a nearest-neighbour
low-order stencil in \citep{Guermond_Nazarov_Popov_2024}.

A different line of work writes the scheme in flux form. Abgrall
\citep[Lemma~4.2]{Abgrall_2018} showed that residual distribution schemes can be
written in flux form: the zero-sum element residuals are turned into edge fluxes
by a small Laplacian solve on each element.
The idea was later used in \citep{Kuzmin_Quezada_2020, Hajduk_etal_2020a,
  Hajduk_etal_2020b, Hajduk_2021} to build subcell fluxes from zero-sum element
contributions; in \citep{Kuzmin_Quezada_2020}, by an element
Laplacian weighted with the $\polP_1$ subcell masses, with one Bernstein mass matrix
for both the low- and high-order schemes. Alternatively, the zero-sum
element contributions can be limited directly, without edge fluxes, which
relates convex limiting to residual distribution, see the review
\citep{Kuzmin_2026}. For discontinuous Galerkin (DG)
methods, we refer to \citep{Vilar_Abgrall_2024}, for an a posteriori local subcell correction
of DG schemes
on unstructured meshes, and to \citep{Pazner_2021}, for sparse IDP DG methods with subcell
convex limiting.

Conservative redistribution (or repair) techniques presented in
\citep{Kucharik_Shashkov_Wendroff_2003, Shashkov_Wendroff_2004,
  Loubere_Staley_Wendroff_2006, Despres_Loubere_2006} are an alternative to FCT. Instead of blending a low- and a high-order scheme,
these methods move mass between nodal states after each update.
The recent works \citep{Guermond_Wang_2025, Guermond_2026}
propose a cell-based Jacobi redistribution on the $\polP_1$ subdivision
of the high-order cells for radiation transport and the compressible Euler
equations, respectively, with continuous finite elements. However, this approach requires
nonnegative lumped masses, which restricts it to $\polP_1$, $\polP_2$ and $\polP_3$ Lagrange
elements on triangles. Moreover, the local redistribution is not guaranteed to
reach the bounds in a fixed number of iterations, and a global redistribution
step is used as a fallback \citep{Guermond_2026}. We refer to \citep{Liu_Cheng_Shu_2017}
and \citep{Griffin_Shu_2026} for a sweeping method that makes one
forward and one backward pass and only changes the values
that violate the bounds.

A graph-Poisson flux reconstruction
for continuous finite elements was introduced in our recent work \citep{Nazarov_2026} for
scalar conservation laws.
Instead of splitting element residuals, it solves one global graph-Poisson
problem whose source is the difference between the high- and low-order updates,
each multiplied by its own lumped mass. This source has zero sum, and the fluxes
are the differences of the solution across the edges.
Fluxes given by differences of nodal potentials across the edges also
appear in the optimization-based flux correction of
\citep{Ruppenthal_Kuzmin_2023}, where the potentials are the control variables
of an inequality-constrained optimization problem; here, the potentials solve a
graph-Poisson problem, globally or element by element, and the resulting fluxes
are limited in closed form.
Together with a low-order scheme on a fine mesh, this construction addresses
both obstacles above. The low-order scheme uses
$\polP_1$ elements on the fine mesh built from the nodes of the high-order
$\polP_k$ mesh. Its stencil is the nearest-neighbour stencil of the fine mesh,
so its viscosity and its time-step bound scale with the local node spacing, as
for $\polP_1$ elements. A first-order viscosity on this fine $\polP_1$ mesh was
introduced in \citep{Dao_Nazarov_2024}, where it bounds the residual viscosity
of a high-order method for magnetohydrodynamics, and a sparsified low-order
operator on the subcell grid is used for DG methods in \citep{Pazner_2021}. The
difference between the high-order mass and the positive $\polP_1$ mass is
carried by the graph-Poisson fluxes, so the low-order fluxes do not have to be
modified as in \citep{Guermond_Nazarov_Popov_2024}, and the invariant domain
property holds whatever the sign of the high-order lumped masses.
Moreover, the graph-Poisson problem only needs a source with zero sum, so
the high-order method only has to be conservative: its stabilization, for
example a continuous interior penalty, does not have to be written in edge form.

This paper aims to address the two obstacles discussed above for the compressible
Euler equations. We list the main contributions of this paper as follows.

\paragraph{Extension to the Euler equations.} We extend the graph-Poisson convex limiting of \citep{Nazarov_2026}
from scalar conservation laws to the compressible Euler equations. The new limiter enforces
local bounds on the density, positivity of the internal energy,
and a local minimum principle on the specific entropy by a line search in a
convex set of admissible states, see Section~\ref{Sec:convex_limiting}.
Boundary data are imposed at prescribed nodes, which are excluded from the
limiting; in the global graph-Poisson problem they enter as a homogeneous
Dirichlet condition, see Section~\ref{Sec:global}.

\paragraph{Element-local fluxes.} We construct element-local graph-Poisson fluxes, in the spirit of
\citep{Abgrall_2018} and \citep{Kuzmin_Quezada_2020}, but for the difference of
two different masses: the fine $\polP_1$ and the coarse $\polP_k$ masses.
Each element contributes a zero-sum vector built from the high-order
mass matrix and the differences of the high- and low-order convection and
stabilization terms, and a small graph-Poisson problem, whose matrix
depends only on the polynomial degree, is solved on each element, see
Section~\ref{Sec:element_local}. No global linear system is solved
for the flux reconstruction.

\paragraph{Zero and negative lumped masses.}
The method is IDP also for Lagrange elements whose lumped
masses are zero or negative at some nodes, such as $\polP_2$ and $\polP_4$
on triangles; the methods of \citep{Guermond_Nazarov_Popov_2024} and
  \citep{Guermond_2026} assume nonnegative lumped masses, which excludes
  $\polP_4$. For such elements, the
  limited update is written with the positive fine $\polP_1$ masses, and an
  element-local exchange restores the high-order mass within the bounds;
  a few more rounds on growing patches and, if needed, a final global round as
  in \citep{Guermond_2026} complete it,
  see Sections~\ref{Sec:HO_LO_relation} and~\ref{Sec:exchange}.

  \paragraph{Numerical validation.}
We test the method with Lagrange elements $\polP_1$--$\polP_4$ on
unstructured meshes, see Section~\ref{Sec:numerics}. Near-optimal convergence
  rates are obtained for smooth problems, and robustness is tested on problems
  with shocks and contact discontinuities, including a Mach 800
  astrophysical jet. 
  The element-local fluxes give the same errors as the global ones for the
  isentropic vortex and are $1.2$--$2.1$ times faster on meshes with about
  $5\cdot10^5$ nodes. For a vortex that leaves the domain through Dirichlet
  boundaries, they keep the convergence rate, whereas the global fluxes lose it
  on fine meshes.

\vspace{0.1in}
  The paper is organized as follows. In Section~\ref{Sec:prelim}
we present the compressible Euler equations, the meshes
and the finite element spaces;
in Section~\ref{Sec:FEM} we present the low- and high-order finite element
approximations of the Euler system; in Section~\ref{Sec:limiting}
we present the global and element-local
graph-Poisson flux reconstructions;
in Section~\ref{Sec:convex_limiting} we present the convex limiting
algorithm; and
in Section~\ref{Sec:numerics} we validate the method on several
challenging problems for the compressible Euler equations.
Section~\ref{Sec:conclusion}
concludes the paper.

\section{Preliminaries}
\label{Sec:prelim}
In this section, we present the model problem and the finite element
meshes and spaces used in the numerical approximation.

\subsection{Compressible Euler equations}
\label{Sec:PDE}
Let $d\in\{1,2,3\}$ be the space dimension and let $\Omega\subset\polR^d$ be
an open bounded polyhedral domain. We seek the conservative variables
$
  \bu := (\rho, \bm, E)^\top \in \polR^{d+2},
$
where $\rho$ denotes the density, $\bm = \rho \bv$ the momentum,
$\bv = (v_1,\ldots,v_d)^\top$ the velocity and $E$ the total energy, such that
\begin{equation}\label{eq:euler}
  \begin{cases}
    \p_t \bu + \DIV \polf(\bu) = 0, & (\bx, t) \in \Omega\times(0,\widehat t\ ], \\
    \bu(\bx, 0) = \bu_0(\bx), & \bx\in\Omega,
  \end{cases}
\end{equation}
with appropriate boundary conditions. Here $\widehat t\ > 0$ is the final time and  the flux is
\begin{equation}\label{eq:flux}
  \polf(\bu) :=
  \Big(
  \bm, \
  \bm \otimes \bv + p\polI, \
  \bv (E+p)
  \Big)^\top \in\polR^{(d+2)\times d},
\end{equation}
and the pressure is given by the ideal gas law
\begin{equation}\label{eq:euler_pressure}
  p(\bu)=(\gamma-1)
  \Big(E-\frac{|\bm|^2}{2\rho}\Big),
  \qquad \gamma>1 .
\end{equation}
For a vector $\bc\in\polR^d$ we write $\polf(\bu)\bc\in\polR^{d+2}$ for the
flux in the direction $\bc$. For $\rho>0$, $p\ge0$, a velocity $\bv$ and a
unit vector $\bn$, the speed of sound and the largest wave speed in the
direction $\bn$ are
\begin{equation}\label{eq:wave_speed}
  \begin{split}
    a(\rho,p):=\sqrt{\gamma p/\rho},
    \quad
    \lambda(\rho,\bv,p;\bn):=|\bv\SCAL\bn|+a(\rho,p),\\
    \quad
    \lambda(\rho,\bv,p):=|\bv|+a(\rho,p)=\max_{|\bn|=1}\lambda(\rho,\bv,p;\bn);
  \end{split}
\end{equation}
$\lambda(\rho,\bv,p;\bn)$ is the spectral radius of the Jacobian of
$\polf(\bu)\bn$. We write $a(\bu)$, $\lambda(\bu;\bn)$ and
$\lambda(\bu)$ for the values at a state $\bu$.

We define the specific internal energy and the specific entropy (up to an
additive constant and a factor $\gamma-1$) by
\[
  \varepsilon(\bu):=\frac{1}{\rho}\Big(E-\frac{|\bm|^2}{2\rho}\Big),
  \qquad
  s(\bu) := \log p(\bu) - \gamma\log\rho .
\]
We use the convex set
\begin{equation}\label{eq:euler_admissible_set}
  \calA
  :=
  \left\{
    \bu=(\rho,\bm,E)^\top\in\mathbb R^{d+2}:\
    \rho>0,\
    \varepsilon(\bu)>0
  \right\},
\end{equation}
which is an invariant domain of \eqref{eq:euler}, and for $s_{\min}\in\polR$
the convex subsets
$\calA_{s_{\min}}:=\{\bu\in\calA:\ s(\bu)\ge s_{\min}\}$.
Equivalently, $\bu\in\calA$ if and only if $\rho>0$ and $p(\bu)>0$.
The mathematical entropy $-\rho s(\bu)$ is a convex function of the
conservative variables on $\calA$
\citep[Thm.~2.1]{Harten_Lax_Levermore_Morojoff_1998}. Since
$s(\bu)\ge s_{\min}$ is equivalent to $\rho s(\bu)-s_{\min}\rho\ge0$ for
$\rho>0$, the sets $\calA_{s_{\min}}$ are convex, that is, $s$ is
quasiconcave on $\calA$; this is the property used by the limiter
\citep{Guermond_Nazarov_Popov_Tomas_2018}.

\subsection{Meshes and finite element spaces}
\label{Sec:spaces}
Let $\calT_h$ be a shape-regular conforming triangulation of $\Omega$ into
simplices $K$, with the affine maps $\bg_K:\widehat K\to K$. For $k\ge1$ we
use the continuous Lagrange space with equispaced nodes
\begin{equation}\label{eq:Xh}
  \calX_h := \{w_h \in \calC^0(\overline \Omega):\ w_h|_K \circ \bg_K \in \polP_k\ \ \forall K\in \calT_h\},
  \qquad
  \bcalW_h := (\calX_h)^{d+2}.
\end{equation}
We denote by $\{\bN_i\}_{i\in\calV}$ the Lagrange nodes, by
$\{\varphi_i\}_{i\in\calV}$ the shape functions, and by $\calI(K)$ the set of
the $n_k:=\binom{k+d}{d}$ nodes of $K$.

\begin{figure}[htbp]
  \centering
  \includegraphics[width=0.8\textwidth]{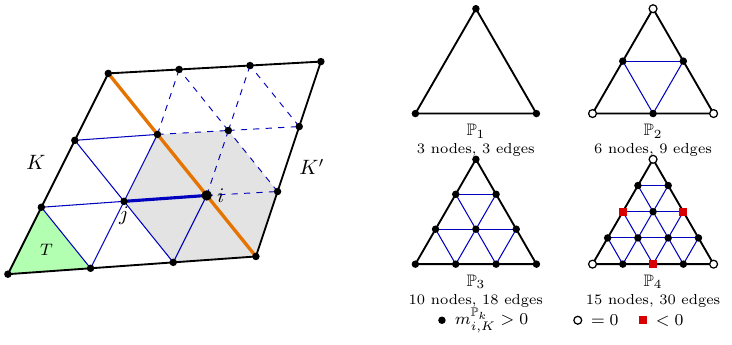}
  \caption{The fine $\polP_1$ mesh. Left: two $\polP_3$ elements $K$ and $K'$
    split into fine simplices, with a fine simplex $T\subset K$ (green). The
    fine edges inside $K$ are solid and those inside $K'$ dashed; the fine
    edges on the common facet (orange) belong to both $\calE(K)$ and
    $\calE(K')$. The support of $\varphi_i^{\polP_1}$ (grey) of a node $i$ on
    the facet consists of fine simplices of both elements, and its other
    vertices form $\calN(i)$; $\{i,j\}\in\calE(K)$ is drawn thick. Right: the
    element graphs $G_K$ for $\polP_1$--$\polP_4$. The marks show the sign of
    the element masses $m_{i,K}^{\polP_k}$ of Section~\ref{Sec:masses}, which
    vanish at the vertices for $\polP_2$ and $\polP_4$ and are negative at the
    edge midpoints for $\polP_4$; the fine masses $m_{i,K}^{\polP_1}$ are
    positive at all nodes.}
  \label{fig:fine_mesh}
\end{figure}

Splitting every $K\in\calT_h$ into $k^d$ simplices whose vertices are the
nodes of $K$ gives the fine $\polP_1$ mesh $\calT_h^{\polP_1}$, with the
$\polP_1$ shape functions $\varphi_i^{\polP_1}$, $i\in\calV$. We write $T$ for
a fine simplex, $T\subset K$ for the fine simplices of a coarse element, and
$\calE$ for the set of edges $\{i,j\}$ of $\calT_h^{\polP_1}$; $\calN(i)$ is
the set of neighbours of $i$ in this graph. Every fine edge lies in one coarse
element (edges inside $K$) or in the elements sharing a coarse facet (edges on
that facet). For a coarse element $K$, $\calE(K)$ denotes the fine edges of
the simplices $T\subset K$; the graph $G_K:=(\calI(K),\calE(K))$ is connected
and has the same structure for all $K$. For $d=2$ it has $n_k=(k+1)(k+2)/2$
nodes and $3k(k+1)/2$ edges, see Figure~\ref{fig:fine_mesh}.

\subsection{Two discrete masses}
\label{Sec:masses}
The consistent mass matrix and its element parts are
\[
  m_{ij}:=\int_\Omega\varphi_j\varphi_i\,\ud\bx=\sum_{K\in\calT_h}m_{ij}^K,
  \qquad
  m_{ij}^K:=\int_K\varphi_j\varphi_i\,\ud\bx,
\]
and the two lumped masses are
\[
  \mk_i:=\sum_{j\in\calV}m_{ij}=\int_\Omega\varphi_i\,\ud\bx,
  \qquad
  \mone_i:=\int_\Omega\varphi_i^{\polP_1}\,\ud\bx .
\]
The fine masses are positive, while $\mk_i$ vanishes at the vertices for
$\polP_2$ and $\polP_4$ and is negative at the edge midpoints for $\polP_4$ on triangles
\citep{Nazarov_2026}. Both split over the elements as $\mone_i=\sum_{K\ni i}m_{i,K}^{\polP_1}$ and $\mk_i=\sum_{K\ni i}m_{i,K}^{\polP_k}$, with the element contributions $m_{i,K}^{\polP_1}:=\int_K\varphi_i^{\polP_1}\,\ud\bx$ and $m_{i,K}^{\polP_k}:=\int_K\varphi_i\,\ud\bx$; moreover, 
$\sum_{i\in\calI(K)}m_{i,K}^{\polP_1}=\sum_{i\in\calI(K)}m_{i,K}^{\polP_k}=|K|$,
and we write
$\Mone(\bsfW):=\sum_{i\in\calV}\mone_i\bsfW_i$, $\Mk(\bsfW):=\sum_{i\in\calV}\mk_i\bsfW_i$.
For a fine simplex $T$ and $i,j\in\calI(T)$ we use the vectors
\[
  \bc_{ij}^T := \int_T \varphi_i^{\polP_1}\GRAD\varphi_j^{\polP_1}\,\ud\bx,
  \qquad
  \bc_{ij}:=\sum_{T\ni i,j}\bc_{ij}^T,
\]
which satisfy $\sum_{j\in\calI(T)}\bc_{ij}^T=\boldsymbol 0$ and
$\bc_{ij}+\bc_{ji}=\int_{\p\Omega}\varphi_i^{\polP_1}\varphi_j^{\polP_1}\bn\,\ud s$;
in particular $\bc_{ij}=-\bc_{ji}$ if $i$ or $j$ is an interior node. We set
$\bn_{ij}:=\bc_{ij}/|\bc_{ij}|$.

\subsection{Boundary conditions and prescribed nodes}
\label{Sec:boundary}
We consider Dirichlet (in particular inflow), supersonic outflow, slip-wall and periodic boundaries.
Boundary conditions are imposed strongly at the nodes, after every
Runge--Kutta stage: Dirichlet values are prescribed, and at a slip-wall node
with unit outward normal $\bn_i$ the normal momentum is removed,
$\bm_i\leftarrow\bm_i-(\bm_i\SCAL\bn_i)\bn_i$, which keeps $\rho_i$, the
tangential momentum and $E_i$, and therefore does not decrease $p_i$
\citep{Guermond_Nazarov_Popov_Tomas_2018}. On curved walls $\bn_i$ is the exact
normal of the physical boundary at $\bN_i$. On periodic boundaries each node is identified with its periodic image, i.e., both copies are a single free node.

We call a node \emph{prescribed} if its state is set by boundary data,
that is, if it lies on an inflow or Dirichlet part of $\p\Omega$, and \emph{free}
otherwise; $\Vd$ and $\Vf$ denote the two sets. Interior nodes, slip-wall nodes and outflow nodes
are free. $\Ef:=\big\{\{i,j\}\in\calE:\ i,j\in\Vf \big\}$ denotes the free edges, and
$\Nf(i):=\calN(i)\cap\Vf$. The limiter acts only on free nodes and free
edges. A prescribed node takes the boundary data after every stage. At an
outflow node nothing is imposed, and the limited state is kept. At a
slip-wall node the limited state is computed as at an interior node, and the
slip projection above is applied to it afterwards. Since the projection keeps
$\rho_i$ and $E_i$ and only lowers the kinetic energy, it does not decrease
$p_i$ or the specific entropy $s_i$, so the bounds enforced by the limiter
still hold after it. For the jet in Section~\ref{Sec:jet}, we also use a
far-field condition and open-boundary nodes whose data are computed after
every stage from the low-order solution of Section~\ref{Sec:FEM}. Both are
treated as prescribed nodes; unlike at Dirichlet nodes, the high-order fluxes
on the edges to them are dropped, so that these edges remain low order, see
Section~\ref{Sec:limited_update}.

\section{Finite element approximation}
\label{Sec:FEM}

\subsection{Low-order scheme}
\label{Sec:first_order}
Let $\bu_h^n=\sum_{i\in\calV}\bsfU_i^n\varphi_i$ be the approximation of $\bu(t_n)$ at the current time $t_n$ and $\tau_n =t_{n+1}-t_n$ the
time step; every stage of the strong stability preserving Runge--Kutta
method used below is a forward Euler step of this form. The low-order method
is the graph-viscosity method of \citep{Guermond_Popov_Fast_Riemann_2016} on
the fine $\polP_1$ graph:
\begin{equation}\label{eq:LO}
  \mone_i\,
  \frac{\bsfU_i^{{\rm L}, n+1}-\bsfU_i^n}{\tau_n}
  =
  -\sum_{j\in\calN(i)\cup\{i\}}\polf(\bsfU_j^n)\,\bc_{ij}
  +
  \sum_{j\in\calN(i)} d_{ij}^{{\rm L}, n}\big(\bsfU_j^n-\bsfU_i^n\big),
  \qquad i\in\calV.
\end{equation}
The first sum
is $-\int_\Omega\DIV(\IoneF)\varphi_i^{\polP_1}\,\ud\bx$, with
$\IoneF:=\sum_{j\in\calV}\polf(\bsfU_j^n)\varphi_j^{\polP_1}$ the fine $\polP_1$
interpolant of the flux, and in the implementation it is assembled fine
simplex by fine simplex from the nodal fluxes. The graph viscosity is
\begin{equation}\label{eq:euler_dijL}
  d_{ij}^{{\rm L}, n}
  :=
  \max\Big\{
  |\bc_{ij}|\,\lambda_{\max}\big(\bn_{ij},\bsfU_i^n,\bsfU_j^n\big),\
  |\bc_{ji}|\,\lambda_{\max}\big(\bn_{ji},\bsfU_j^n,\bsfU_i^n\big)
  \Big\}
  = d_{ji}^{{\rm L},n},
\end{equation}
where $\lambda_{\max}(\bn,\bsfU_L^n,\bsfU_R^n)$ is an upper bound on the maximum
wave speed of the one-dimensional Riemann problem in the direction $\bn$ with
data $\bsfU_L^n$, $\bsfU_R^n$; in particular
$\lambda_{\max}(\bn,\bsfU_L^n,\bsfU_R^n)\ge\max\{\lambda(\bsfU_L^n;\bn),\lambda(\bsfU_R^n;\bn)\}$.
A guaranteed bound is obtained from the
two-rarefaction estimate or from Toro's pressure equation
\citep{Guermond_Popov_Fast_Riemann_2016,Toro_2009}.

Since $\sum_{j\in\calV}\bc_{ij}=\boldsymbol 0$, \eqref{eq:LO} can be written in the
bar-state form
\begin{equation}\label{eq:LO_bar}
  \bsfU_i^{{\rm L}, n+1}
  =
  \Big(1-\frac{2\tau_{n}}{\mone_i}\sum_{j\in\calN(i)}d_{ij}^{{\rm L}, n}\Big)\bsfU_i^n
  +\frac{2\tau_{n}}{\mone_i}\sum_{j\in\calN(i)}d_{ij}^{{\rm L}, n}\,\Ubar_{ij},
  \qquad
  \Ubar_{ij}:=\frac{\bsfU_i^n+\bsfU_j^n}{2}
  -
  \big(\polf(\bsfU_j^n)-\polf(\bsfU_i^n)\big)\frac{\bc_{ij}}{2d_{ij}^{{\rm L}, n}} ,
\end{equation}
for every interior node $i$; at boundary nodes the boundary term
$\polf(\bsfU_i^n)\sum_{j\in\calN(i)\cup\{i\}}\bc_{ij}$ is replaced by the boundary condition.

\begin{proposition}[Low-order method]\label{prop:LO}
  The update \eqref{eq:LO} is conservative, \ie
  \[
    \Mone(\bsfU^{{\rm L}, n+1})=\Mone(\bsfU^n)-\tau_{n}\int_{\p\Omega}\IoneF\bn\,\ud s.
  \]
  If $\bsfU_j^n\in\calA$ for all $j$ and
  \begin{equation}\label{eq:cfl}
    \frac{2\tau_n}{\mone_i}\sum_{j\in\calN(i)}d_{ij}^{{\rm L}, n}\le 1
    \qquad\forall i\in\calV,
  \end{equation}
  then $\bsfU_i^{{\rm L}, n+1}$ is a convex combination of $\bsfU_i^n$ and the bar
  states $\Ubar_{ij}$. With a guaranteed wave-speed bound,
  $\Ubar_{ij}\in\calA$ and $s(\Ubar_{ij})\ge\min\{s(\bsfU_i^n),s(\bsfU_j^n)\}$,
  so that $\bsfU_i^{{\rm L}, n+1}\in\calA$ and
  $s(\UL_{i})\ge\min_{j\in\calN(i)\cup\{i\}}s(\bsfU_j^n)$.
\end{proposition}
\begin{proof}
  Summing \eqref{eq:LO} over $i$, the viscous terms cancel since
  $d_{ij}^{{\rm L}, n}=d_{ji}^{{\rm L}, n}$, and
  $\sum_{i\in \calV}\sum_{j\in\calN(i)\cup\{i\}}\polf(\bsfU_j^n)\bc_{ij}=\int_\Omega\DIV(\IoneF)\,\ud\bx$ by the
  partition of unity. The convex combination is \eqref{eq:LO_bar}; its
  coefficients are nonnegative by \eqref{eq:cfl} and sum to one. The
  properties of the bar states are those of the averaged Riemann solution
  \citep[Sec.~2]{Guermond_Popov_Fast_Riemann_2016}, and $\calA_{s_{\min}}$ is
  convex.
\end{proof}

\subsection{High-order scheme}
\label{Sec:euler_high_order}

The high-order method adapts the residual viscosity (RV) method of
\citep{Nazarov_Hoffman_2013} to the fine graph. The RV method captures
shocks and contact discontinuities, but, like the underlying Galerkin
method, it converges at a suboptimal rate for polynomial spaces of even
degree; see \citep[Sec.~5]{Dao_Nazarov_2022} and the references therein.
To recover the optimal rate for even degrees
\citep[Sec.~4.2]{Ern_Guermond_2013}, we add the continuous interior
penalty (CIP) stabilization of \citep{Douglas_Dupont_1976,Burman_Hansbo_2004}.

\subsection*{Galerkin approximation in strong form.}
Let $\bu_h^n=\sum_{i\in\calV}\bsfU_i^n\varphi_i=:(\rho_h^n,\bm_h^n,E_h^n)^\top$, and let
\[
  \bv_h^n:=\sum_{i\in\calV}\bv_i^n\varphi_i,\quad \bv_i^n:=\frac{\bm_i^n}{\rho_i^n},
  \qquad
  p_h^n:=\sum_{i\in\calV}p(\bsfU_i^n)\varphi_i
\]
be the $\polP_k$ interpolants of the nodal velocities and pressures. Thus
$\bm_i^n=\rho_i^n\bv_i^n$ at the nodes, whereas $\bm_h^n\neq\rho_h^n\bv_h^n$ in
general. We approximate the flux by products of these interpolants,
\begin{equation}\label{eq:Fh}
  \Fh
  :=
  \Big(\bm_h^n,\ \bm_h^n\otimes\bv_h^n+p_h^n\polI,\ (E_h^n+p_h^n)\bv_h^n\Big)^\top .
\end{equation}
The high-order candidate $\bu_h^{{\rm H}, n+1}:=\sum_{i\in\calV}\bsfU_i^{{\rm H}, n+1}\varphi_i$ solves
\begin{equation}\label{eq:HO}
  \sum_{j\in\calV} m_{ij}
  \frac{\bsfU_j^{{\rm H}, n+1} - \bsfU_j^{n}}{\tau_n}
  = \bsfb_i^n
  :=
  -\int_{\Omega} \DIV \Fh\, \varphi_i \,\ud \bx
  +\sum_{j\in\calN(i)}d_{ij}^{{\rm H}, n}\big(\bsfU_j^n-\bsfU_i^n\big)
  -j_h\big(\bu_h^n,\varphi_i\big),
  \qquad i \in \calV,
\end{equation}
with the graph residual viscosity $d_{ij}^{{\rm H}, n}$ of
Section~\ref{Sec:rv} and the continuous interior penalty term $j_h$ of
Section~\ref{Sec:cip}. The convective term is written in strong form: no
boundary integral appears. Summing over $i$, the viscous and penalty terms
cancel, and what remains is the boundary flux
$-\int_{\p\Omega}\Fh\bn\,\ud s$.

For all $k$ we solve the linear system \eqref{eq:HO} exactly: the consistent
mass matrix is block diagonal over the $d+2$ components with the same scalar
block, which does not change in time and is factorized once. Hence
$\tau_n\bsfb_i^n=\sum_{j\in\calV}m_{ij}(\bsfU_j^{{\rm H},n+1}-\bsfU_j^n)$ holds at
every node, up to round-off.

For $k=1$ we replace \eqref{eq:Fh} by its linear interpolant
$\Fh:=\IoneF$ as in the low-order scheme \eqref{eq:LO}.

\subsection*{Continuous interior penalty.}
\label{Sec:cip}
For an interior facet $F$ with diameter $h_F$ and unit normal $\bn_F$, let
$K^+,K^-\in\calT_h$ be the two elements sharing $F$, and let
\[
  \lambda_F^n:=\lambda(\rho_h^n,\bv_h^n,p_h^n;\bn_F)\qquad\text{on }F,
\]
the wave speed \eqref{eq:wave_speed} of the interpolants $\rho_h^n$, $\bv_h^n$,
$p_h^n$ defined before \eqref{eq:Fh}. These are continuous, so the traces from
$K^+$ and $K^-$ agree and $\lambda_F$ needs no average. We define,
componentwise,
\begin{equation}\label{eq:cip}
  j_h(\bu_h^n, \bw)
  :=
  \sum_{F\in\calF_h^{\rm int}}
  j_F(\bu_h^n,\bw),
  \qquad
  j_F(\bu_h^n,\bw):=
  \gamma_{\rm cip} h_F^2
  \int_F
  \lambda_F^{n}\,
  \jump{\GRAD \bu_h^n\bn_F}
  \SCAL
  \jump{\GRAD \bw \bn_F}
  \,\ud s ,
\end{equation}
following \citep{Douglas_Dupont_1976,Burman_Hansbo_2004}. The facet term
$j_F(\bu_h^n,\varphi_i)$ vanishes unless $i\in\calI(K^+)\cup\calI(K^-)$, the
nodes of the two elements sharing $F$, and $\sum_{i\in\calV}j_F(\bu_h^n,\varphi_i)=0$
since the jump of the gradient of a constant vanishes. Following
\citep[Eq.~(4.6)]{Guermond_Wang_2025}, we set
$\gamma_{\rm cip}:=c_{\rm cip}\,\varpi$ with $\varpi:=d^2/(1+k)^4$, and we use
$c_{\rm cip}=0.05$ unless stated otherwise, for all $k$. 

\subsection*{Graph residual viscosity.}
\label{Sec:rv}
The high-order viscosity is a fraction of the low-order one on every edge,
\begin{equation}\label{eq:euler_dijH}
  d_{ij}^{{\rm H}, n}:=d_{ij}^{{\rm L}, n}\max\{r_i^n,r_j^n\},
  \qquad
  r_i^n:=\frac{\mu_i^{{\rm H}, n}}{\mu_i^{{\rm L}, n}}\in[0,1],
\end{equation}
with the first-order and residual-based nodal coefficients
\begin{equation}\label{eq:mu}
  \mu_i^{{\rm L}, n}:=C_{\max}\,\lambda_i^n\,h_i,
  \qquad
  \mu_i^{{\rm H}, n}:=\min\big\{\mu_i^{{\rm L}, n},\ C_{\rm R}\,h_i^2\,R^n_i\big\}.
\end{equation}
Here $C_{\textrm{R}}>0$ is a user-defined parameter,
$\lambda_i^n$ and $h_i$ are patch quantities on the fine graph,
\begin{equation}\label{eq:lambda_h}
  \lambda_i^n:=\max_{j\in\{i\}\cup\calN(i)}\lambda(\bsfU_j^n),\qquad
  h_i:=\Big(\max_{j\in\{i\}\cup\calN(i)}\ \max_{T\ni j}\|\GRAD\varphi_j^{\polP_1}|_T\|\Big)^{-1},
\end{equation}
and $R_i^n:=\max_{j\in\{i\}\cup\calN(i)}\max\{0,\widehat{\mathcal R}^n_j\}$
is the patch maximum of the lumped $L^2$-projection $\widehat{\mathcal R}^n_h \in \calX_h$ of
the normalized residual
\begin{equation}\label{eq:euler_residual_norm}
  \mathcal R^n
  :=
  \max\left\{
  \frac{|D_t^n\bm_h+\DIV(\bm^n_h\otimes\bv^n_h)+\GRAD p^n_h|}{\operatorname{osc}(\bm^n_h)},\
  \frac{|D_t^n\rho_h+\DIV(\rho^n_h\bv^n_h)|}{\operatorname{osc}(\rho^n_h)},\
  \frac{|D_t^nE_h+\DIV(\bv_h(E_h^n+p_h^n))|}{\operatorname{osc}(E_h^n)}
  \right\},
\end{equation}
normalized by the spread of the nodal values around their mean,
\[
  \operatorname{osc}(f_h^n):=\max\Big\{\max_{i\in\calV}\big|\mathsf f^n_i-\bar{\mathsf f}^n\big|,\ 10^{-12}\Big\},
  \qquad
  \bar{\mathsf f}^n:=\frac{1}{|\calV|}\sum_{i\in\calV}\mathsf f^n_i ,
\]
with the Euclidean norm and the componentwise mean for the momentum, so that
$\mathcal R^n$ is invariant under rotations of the coordinates. It is evaluated at the quadrature points (degree $3k-1$) with the $\polP_k$
interpolants $\bv_h$, $p_h$ of the nodal velocity and pressure and, for
$D_t$, the variable-step backward difference through
$\bu^{n}_h,\dots,\bu^{n-q}_h$ of order $q=\min\{n,4\}$:
$\widehat{\mathcal R}^n_i:=\frac{1}{\mk_i}\int_\Omega\mathcal R^n \varphi_i\,\ud \bx$ where $\mk_i>0$, and
the mean of $\widehat{\mathcal R}^n$ over the fine-graph neighbours with $\mk>0$
where $\mk_i\le0$ (the vertices for $\polP_2$, the vertices and edge midpoints for $\polP_4$). The residual is computed once per time step. In smooth
regions $\mathcal R^n=\mathcal O(h^k)$, since it contains first derivatives of $\bu_h$, hence $r_i^n=\mathcal O(h^{k+1})$ and the high-order viscosity perturbs the scheme by $\mathcal O(h^{k+1})$. Near shocks
$r_i^n=1$ and $d_{ij}^{{\rm H}, n}=d_{ij}^{{\rm L}, n}$.

\begin{remark}[Start-up of the time derivative]
  \label{rem:bdf_startup}
  In the start-up phase the states $\bu^{-1}_h,\bu^{-2}_h,\bu^{-3}_h$
  are not available. In the first three time steps we therefore replace the
  fourth-order formula by the backward difference formulas of order one, two
  and three, \ie we use the order $q=\min\{n+1,4\}$ at $t^{n+1}$.
\end{remark}

\begin{proposition}[High-order method]\label{prop:HO}
  $\Mk(\bsfU^{{\rm H}, n+1})=\Mk(\bsfU^n)-\tau_{n}\int_{\p\Omega}\Fh\bn\,\ud s$.
\end{proposition}
\begin{proof}
  Sum \eqref{eq:HO} over $i$; the viscous and CIP terms sum to zero.
\end{proof}

\section{Graph-Poisson flux reconstruction}
\label{Sec:limiting}
By construction, the low- and high-order schemes of Section~\ref{Sec:FEM} conserve the masses $\Mone$ and $\Mk$, respectively; see Propositions~\ref{prop:LO} and~\ref{prop:HO}. As in \citep{Nazarov_2026}, we compare the two updates, each in its own mass, componentwise:
\begin{equation}\label{eq:Ai}
  \bsfA_i^{n+1}
  :=
  \mk_i\big(\UH_i-\bsfU_i^n\big)
  -
  \mone_i\big(\UL_i-\bsfU_i^n\big),
  \qquad i\in\calV .
\end{equation}

\begin{lemma}\label{lemma:Ai}
  Assume periodic boundary conditions, or $\Fh\bn =\IoneF\bn = 0$ on $\p\Omega$. Then
  $
  \sum_{i\in\calV} \bsfA_i^{n+1} = 0.
  $
\end{lemma}
\begin{proof}
  See \citep[Lemma~5.1]{Nazarov_2026}.
\end{proof}

We want to write $\bsfA_i^{n+1}$ as a sum of antisymmetric edge fluxes $\bsfA_{ij}^{n+1}=-\bsfA_{ji}^{n+1}$,
$\{i,j\}\in\calE$, with
\begin{equation}\label{eq:Aij}
  \bsfA_i^{n+1} = \sum_{j\in\calN(i)}\bsfA_{ij}^{n+1}.
\end{equation}

There are many ways of obtaining such edge fluxes. We take the fluxes of minimal $\ell^2$-norm,
$\bsfA_{ij}=\bsfQ_i-\bsfQ_j$, where $\bsfQ$ solves the graph-Poisson problem $L\bsfQ=\bsfA$, with
the graph Laplacian $(L\bsfQ)_i=\sum_{j\in\calN(i)}(\bsfQ_i-\bsfQ_j)$ of the
fine graph applied componentwise, see \citep[Sec.~5.1]{Nazarov_2026}.

\subsection{Relation between the two updates}\label{Sec:HO_LO_relation}
Before we construct the limiter, we write the high-order update in terms of the low-order one. We distinguish two cases: if all high-order lumped masses are positive, we use the \emph{$\polP_k$-form}; if some of them are zero or negative, we use the \emph{fine $\polP_1$-form}. 

\subsubsection*{$\polP_k$-form.} 
Here $\mk_i>0$ for all $i\in\calV$, and we write \eqref{eq:Ai} as
\[
  \UH_i = \bsfU^n_i + \frac{\mone_i}{\mk_i} \Big(\UL_i - \bsfU^n_i \Big) +
  \frac{1}{\mk_i} \bsfA_i^{n+1}.
\]

Upon setting
$
\ULhat_i:= \Un_i + \frac{\mone_i}{\mk_i} \Big(\UL_i - \Un_i \Big), 
$
and using \eqref{eq:Aij} we obtain
\begin{equation}\label{eq:HO_LO_pk}
  \UH_i = \ULhat_i  +  \frac{1}{\mk_i} \sum_{j\in\calN(i)}\bsfA_{ij}^{n+1}.
\end{equation}

The base state $\ULhat_i$ is the low-order update with the time-step $\tau_n \frac{\mone_i}{\mk_i}$, and it satisfies Proposition \ref{prop:LO} if the CFL condition \eqref{eq:cfl} holds with $\tau_n$ replaced by $\tau_n \frac{\mone_i}{\mk_i}$.

\subsubsection*{Fine $\polP_1$-form.}
Here $\mk_i\le0$ for some $i\in\calV$, so we cannot divide \eqref{eq:Ai} by $\mk_i$. Instead, we divide it by $\mone_i$, add $\UH_i$ to both sides, and rearrange:
\begin{equation}\label{eq:HO_LO_p1}
  \UH_i = \UL_i  +  \frac{1}{\mone_i} \sum_{j\in\calN(i)}\bsfA_{ij}^{n+1} +\bsfs_i^{n+1},
\end{equation}
where
\[
  \bsfs_i^{n+1}:=\frac{\mone_i-\mk_i}{\mone_i}\big(\bsfU_i^{{\rm H},n+1}-\bsfU_i^n\big),
\]
is a nodal lumping source.

\subsection{Global graph-Poisson fluxes}
\label{Sec:global}
Lemma~\ref{lemma:Ai} shows that $\bsfA$ has sum zero for periodic and no-flux boundary conditions, and therefore the graph-Poisson problem $L \bsfQ = \bsfA$ has a solution, unique up to an additive constant that does not change the fluxes. Since this problem couples all nodes, we call these fluxes \emph{global}. 

Recall that the limiter of Section~\ref{Sec:convex_limiting} acts only on the free nodes, whereas the prescribed nodes take the boundary data. Below we show how the graph-Poisson problem is constructed in this case.

\paragraph{Graph Laplacian.} We use the graph Laplacian of the whole fine
graph, as in \citep{Nazarov_2026},
\begin{equation}\label{eq:L_graph}
  (L\bsfQ)_i:=\sum_{j\in\calN(i)}\big(\bsfQ_i-\bsfQ_j\big),\qquad i\in\calV,
\end{equation}
applied to every conserved component. The neighbours of a free node are
either free or prescribed, so at a free node the row splits into
\begin{equation}\label{eq:L_split}
  (L\bsfQ)_i
  =\sum_{j\in\Nf(i)}\big(\bsfQ_i-\bsfQ_j\big)
  +\sum_{b\in\calN(i)\cap\Vd}\big(\bsfQ_i-\bsfQ_b\big),
  \qquad i\in\Vf .
\end{equation}
The first sum will be the flux to the free neighbours, and the second the
flux to the prescribed neighbours. The graph is connected in all our tests;
otherwise everything below applies to each connected component.

\paragraph{Prescribed nodes.}
A prescribed node takes the boundary data after every stage, whatever the scheme has computed there. Anything we send to it is overwritten by the data, so a prescribed node can absorb any flux.

Since the limiter changes only the free nodes, the fluxes have to reproduce
$\bsfA_i^{n+1}$ only there: if $\sum_{j\in\calN(i)}\bsfA_{ij}^{n+1}=\bsfA_i^{n+1}$ at every
free node, then \eqref{eq:HO_LO} gives back $\UH_i$ when nothing is limited.
Fluxes between free nodes alone cannot do this in general. A flux on an edge
between two free nodes adds $\bsfA_{ij}^{n+1}$ at one end and
$-\bsfA_{ij}^{n+1}$ at the other, so the sum of all these fluxes over the free
nodes vanishes, whereas $\sum_{i\in\Vf}\bsfA_i^{n+1}$ does not vanish in
general. The difference has to be carried by something else: by fluxes to the
prescribed nodes, or by sources at the boundary nodes.
Summing
\eqref{eq:Ai} over all nodes and using Propositions~\ref{prop:LO}
and~\ref{prop:HO} for the updates before the boundary conditions are applied
gives, by the partition of unity,
\[
  \sum_{i\in\calV}\bsfA_i^{n+1}=-\tau_n\sum_{i\in\calV}\beta_i^n,
  \qquad
  \beta_i^n:=\int_{\p\Omega}\big(\Fh-\IoneF\big)\bn\,\varphi_i^{\polP_1}\,\ud s .
\]
Here $\tau_n\beta_i^n$ is the difference between the high- and the low-order
flux through the boundary facets at node $i$; it vanishes for $k=1$ and away
from the boundary, and on a straight slip wall its density and energy
components vanish except at the corners.

\paragraph{Boundary terms.}
We treat the boundary terms at the boundary nodes where they arise. Besides $\tau_n\beta_i^n$ there is a second boundary term from the low- and high-order updates: $\UL$ and $\UH$ are computed first, and the nodal boundary conditions are applied to both before they are handed to the limiter,
so that at a slip-wall node $\bsfA_i^{n+1}$ also contains
\[
  \bsfp_i^{n+1}:=\mk_i\,\Delta\UH_i-\mone_i\,\Delta\UL_i ,
\]
where $\Delta$ denotes the change made by the slip projection; it has only a
normal momentum component. At a free node we collect both in the boundary
term
\[
  \bsfB_i^{n+1}:=\tau_n\beta_i^n-\bsfp_i^{n+1},
\]
which vanishes away from the boundary, move it out of the graph-Poisson
problem, and apply $-\bsfB_i^{n+1}$ at node $i$ as a limited nodal source. At
a slip-wall node the density and energy components of $-\bsfB_i^{n+1}$ are not
applied: they are the high-order flux of mass and energy through the wall at
corners and on curved walls, which is physically zero.

The remaining sum $\sum_{i\in\Vf}(\bsfA_i^{n+1}+\bsfB_i^{n+1})$ is what the two
updates exchange with the prescribed nodes, including the difference of the
boundary fluxes at the inflow boundary. If $\Vd\neq\emptyset$, it may leave the
free nodes through the prescribed nodes, since whatever arrives there is
overwritten by the data. If $\Vd=\emptyset$, all nodes are free, and the sum
vanishes, since the mass system is solved exactly. In both cases the
right-hand side of the graph-Poisson problem is
\begin{equation}\label{eq:A_compatible}
  \widehat\bsfA_i^{n+1}:=\bsfA_i^{n+1}+\bsfB_i^{n+1},
  \qquad i\in\Vf .
\end{equation}
Removing the sum as a uniform shift instead, as in \citep{Nazarov_2026}, would
change every state by $\sum_{i\in\Vf}\bsfA_i^{n+1}/\sum_{j\in\Vf}\mu_j$, with $\mu_j$ of \eqref{eq:HO_LO}, which for $k\ge2$
costs about one order of accuracy.

We take the fluxes of smallest $\ell^2$-norm, on the edges with at least one
free end point, that reproduce $\widehat\bsfA$ at the free nodes. With a
Lagrange multiplier $\bsfQ_i$ for the balance at each free node, the
stationarity condition for the flux on an edge $\{i,j\}$ between free nodes,
which enters the balances at $i$ and at $j$ with opposite signs, is
$\bsfA_{ij}=\bsfQ_i-\bsfQ_j$, as for $L\bsfQ=\bsfA$ \citep{Nazarov_2026}. The
flux on an edge $\{i,b\}$ with $b\in\Vd$ enters only the balance at $i$,
since a prescribed node has no balance to satisfy, and the stationarity
condition gives $\bsfQ_i$, that is, $\bsfQ_i-\bsfQ_b$ with
$\bsfQ_b=\boldsymbol 0$.
By \eqref{eq:L_split}, $\bsfQ$ solves the
graph-Poisson problem with a homogeneous Dirichlet condition at the prescribed
nodes,
\begin{equation}\label{eq:L_dirichlet}
  (L\bsfQ)_i=\widehat\bsfA_i,\quad i\in\Vf,
  \qquad
  \bsfQ_b=\boldsymbol 0,\quad b\in\Vd .
\end{equation}
If $\Vd=\emptyset$, all nodes are free, and \eqref{eq:L_dirichlet} is
$L\bsfQ=\widehat\bsfA$ on the whole graph. It is solvable since
$\sum_{i\in\calV}\widehat\bsfA_i=\boldsymbol 0$, and it determines $\bsfQ$ up
to an additive constant, which does not change the fluxes
$\bsfQ_i-\bsfQ_j$; we fix it by setting $\bsfQ=\boldsymbol 0$ at one node.

\paragraph{Solution of the graph-Poisson problem.}
The matrix of \eqref{eq:L_dirichlet} is symmetric and diagonally dominant, and
strictly so in the rows of free nodes with a prescribed neighbour; hence it is
positive definite if every connected component of the free graph contains
such a node, which is the case in all our tests. The system can be solved fast with multigrid
preconditioners.

\subsection{Element-local fluxes}
\label{Sec:element_local}

The global fluxes of Section~\ref{Sec:global} require a graph-Poisson solve on the whole fine graph in every Runge--Kutta stage, which can be time consuming.

In this section, we construct the fluxes element by element, in the spirit of the conservative residual distribution schemes of \citep{Abgrall_2018}, the local subcell flux reconstruction of  \citep{Kuzmin_Quezada_2020} and \citep{Hajduk_2021}, and the a posteriori subcell correction scheme of \citep{Vilar_Abgrall_2024}.

\subsection*{Splitting $\bsfA$ into element vectors.}
Let $I_h^{\polP_1}\Fh$ denote the fine-$\polP_1$ interpolant of the function $\Fh$ and $\IoneF$ be the fine-$\polP_1$ interpolant of the function $\polf$ as before, \ie
$\IoneF=\sum_{j\in\calV}\polf(\bsfU_j^n)\varphi_j^{\polP_1}$.
Let $\dH:=\bsfU^{{\rm H}, n+1}-\bsfU^n$ denote the increment computed from \eqref{eq:HO} before the
boundary conditions are applied, $\sum_{j\in\calV}m_{ij}\dH_j=\tau_n\bsfb_i^n$.

For every coarse element $K\in\calT_h$ and $i\in \calI(K)$ we define the element vector
\begin{equation}\label{eq:aK}
  \begin{aligned}
    \bsfa_i^{K, n+1}
    :={}&
    \sum_{j\in\calI(K)}m_{ij}^K\big(\dH_i-\dH_j\big)
    &&\text{(mass)}
    \\
    &-\tau_{n}\int_K\DIV\Fh\,\varphi_i\,\ud\bx
    +\tau_{n}\int_K\DIV\big(\IoneF\big)\,\varphi_i^{\polP_1}\,\ud\bx
    &&\text{(high- and low-order convection)}
    \\
    &+\tau_{n}\int_{\p K}\big(\Fh-\IoneF\big)\bn_K\,\varphi_i^{\polP_1}\,\ud s
    &&\text{(facet terms)}
    \\
    &-\tau_{n}\, \bsfj_i^{K,n}
    &&\text{(CIP)},
  \end{aligned}
\end{equation}
where $\bn_K$ is the outward unit normal vector of $K$. The low-order convection is evaluated on the fine simplices $T\in\calT_h^{\polP_1}$ with $T\subset K$,
\[
  \int_K\DIV\big(\IoneF\big)\varphi_i^{\polP_1}\,\ud\bx
  =\sum_{{T\in\calT_h^{\polP_1}:\,T\subset K,\,T\ni i}}\ \sum_{j\in\calI(T)}\big(\polf(\bsfU_j^n)-\polf(\bsfU_i^n)\big)\bc_{ij}^T .
\]

\paragraph{Splitting of the CIP term.} The element vectors must have zero sum
over their element, but the CIP term couples the two elements of an interior
facet, so its facet contributions have to be split between them. Write
$j_h=\sum_F j_F$ as a sum over the interior facets, and for an interior facet
$F$ shared by $K^+$ and $K^-$ let
\[
  \bsfc_i^{F,n}:=j_F(\bu_h^n,\varphi_i),\qquad i\in\calV .
\]
$\bsfc_i^{F,n}$ vanishes unless $i\in\calI(K^+)\cup\calI(K^-)$, and
$\sum_i\bsfc_i^{F,n}=j_F(\bu_h^n,1)=\boldsymbol 0$, since the normal derivative of
a constant has no jump. Let $S_F:=\calI(K^+)\cap\calI(K^-)$ be the nodes on $F$,
$\calI_F^\pm:=\calI(K^\pm)\setminus S_F$ the remaining nodes of $K^\pm$, and
$\bsfc^{+,n}:=\sum_{i\in\calI_F^+}\bsfc_i^{F,n}$. We split $\bsfc^{F,n}$ into
\begin{equation}\label{eq:cip_split}
  \bsfc_i^{F,K^+, n}:=
  \begin{cases}
    \bsfc_i^{F,n}, & i\in\calI_F^+,\\
    -\bsfc^{+,n}/|S_F|, & i\in S_F,
  \end{cases}
  \qquad
  \bsfc_i^{F,K^-,n}:=
  \begin{cases}
    \bsfc_i^{F,n}, & i\in\calI_F^-,\\
    \bsfc_i^{F,n}+\bsfc^{+,n}/|S_F|, & i\in S_F ,
  \end{cases}
\end{equation}
and zero at all other nodes: $K^+$ keeps its entries off $F$ and balances them
at the nodes of $F$, and $K^-$ takes the rest. Each part has zero sum over its
element (for $K^-$ because $\sum_i\bsfc_i^{F,n}=\boldsymbol 0$), and the two parts
add up to $\bsfc^{F,n}$. The CIP part of the element vector of $K$ is
$\bsfj_i^{K,n}:=\sum_{F}\bsfc_i^{F,K,n}$, the sum over the interior facets of $K$.
Which element is $K^+$ does not matter for these properties, but the choice
must not depend on the numbering of the elements, so that the fluxes do not
depend on the mesh partition in parallel; we take as $K^+$ the element whose
vertex opposite to $F$ comes first in the lexicographic order of the
coordinates $(x,y)$.

\begin{lemma}[Zero sum]\label{lem:zero_sum_K}
  For every $K\in\calT_h$, $\sum_{i\in\calI(K)}\bsfa_i^{K, n+1}=\boldsymbol 0$.
\end{lemma}
\begin{proof}
  The mass part vanishes by the symmetry of $m_{ij}^K$. Since
  $\sum_{i\in\calI(K)}\varphi_i=\sum_{i\in\calI(K)}\varphi_i^{\polP_1}=1$ on
  $K$, the divergence theorem gives
  $-\tau_{n}\int_{\p K}\Fh\bn_K\,\ud s$ for the high-order convection,
  $+\tau_{n}\int_{\p K}\IoneF\bn_K\,\ud s$ for the low-order convection, and
  $+\tau_{n}\int_{\p K}(\Fh-\IoneF)\bn_K\,\ud s$ for the facet terms. The three
  cancel. The CIP parts have zero sum by \eqref{eq:cip_split}.
\end{proof}

The facet terms are what makes the element vectors compatible. They are
small for smooth solutions: on a facet, $\Fh-I_h^{\polP_1}\Fh$ is an
interpolation error, and $\IoneF-I_h^{\polP_1}\Fh$ vanishes at the nodes.

\begin{lemma}[Exact splitting]\label{lem:split}
  For every free node $i\in\Vf$,
  \begin{equation}\label{eq:split}
    \widehat\bsfA_i^{n+1}
    =\sum_{K\ni i}\bsfa_i^{K, n+1}
    +\tau_n\sum_{j\in\calN(i)}\big(d_{ij}^{{\rm L},n}-d_{ij}^{{\rm H},n}\big)
    \big(\bsfU_i^n-\bsfU_j^n\big),
  \end{equation}
  with $\widehat\bsfA_i^{n+1}$ from \eqref{eq:A_compatible}.
\end{lemma}
\begin{proof}
  Let $\Delta\UH_i$ and $\Delta\UL_i$ be the changes made by the boundary
  conditions, so that $\UH_i-\bsfU_i^n=\dH_i+\Delta\UH_i$. At a free node they
  vanish except at slip-wall nodes, where
  $\mk_i\Delta\UH_i-\mone_i\Delta\UL_i=\bsfp_i^{n+1}$. Hence
  \[
    \bsfA_i^{n+1}
    =\mk_i\dH_i
    -\mone_i\big(\UL_i-\Delta\UL_i-\bsfU_i^n\big)
    +\bsfp_i^{n+1}.
  \]

  \emph{High-order part.} Since $\mk_i=\sum_{j\in\calV}m_{ij}$ and
  $m_{ij}=\sum_{K\in\calT_h}m_{ij}^K$,
  \[
    \mk_i\dH_i
    =\sum_{j\in\calV}m_{ij}\dH_j+\sum_{j\in\calV}m_{ij}\big(\dH_i-\dH_j\big)
    =\tau_n\bsfb_i^n
    +\sum_{K\ni i}\sum_{j\in\calI(K)}m_{ij}^K\big(\dH_i-\dH_j\big),
  \]
  and by \eqref{eq:HO}
  \[
    \tau_n\bsfb_i^n
    =-\tau_n\sum_{K\ni i}\int_K\DIV\Fh\,\varphi_i\,\ud\bx
    -\tau_n\sum_{j\in\calN(i)}d_{ij}^{{\rm H},n}\big(\bsfU_i^n-\bsfU_j^n\big)
    -\tau_n\sum_{K\ni i}\bsfj_i^{K,n},
  \]
  since the two parts in \eqref{eq:cip_split} add up to $\bsfc^{F,n}$, so that
  $\sum_{K\ni i}\bsfj_i^{K,n}=j_h(\bu_h^n,\varphi_i)$.

  \emph{Low-order part.} $\UL_i-\Delta\UL_i$ is the update \eqref{eq:LO}.
  Its first sum is assembled over the fine simplices, and each of them lies in
  one element $K$, so
  \[
    -\mone_i\big(\UL_i-\Delta\UL_i-\bsfU_i^n\big)
    =\tau_n\sum_{K\ni i}\int_K\DIV\big(\IoneF\big)\varphi_i^{\polP_1}\,\ud\bx
    +\tau_n\sum_{j\in\calN(i)}d_{ij}^{{\rm L},n}\big(\bsfU_i^n-\bsfU_j^n\big).
  \]

  \emph{Facet terms.} They are not part of $\bsfA_i^{n+1}$; we add and
  subtract them. On an interior facet the terms of the two elements cancel,
  since $\Fh$, $\IoneF$ and $\varphi_i^{\polP_1}$ are continuous and
  $\bn_{K^+}=-\bn_{K^-}$. Only the boundary facets remain:
  \[
    \sum_{K\ni i}\tau_n\int_{\p K}\big(\Fh-\IoneF\big)\bn_K\,
    \varphi_i^{\polP_1}\,\ud s=\tau_n\beta_i^n .
  \]

  Collecting the three parts,
  \[
    \bsfA_i^{n+1}
    =\sum_{K\ni i}\bsfa_i^{K, n+1}
    +\tau_n\sum_{j\in\calN(i)}\big(d_{ij}^{{\rm L},n}-d_{ij}^{{\rm H},n}\big)
    \big(\bsfU_i^n-\bsfU_j^n\big)
    -\tau_n\beta_i^n+\bsfp_i^{n+1},
  \]
  which is \eqref{eq:split}, since
  $\bsfB_i^{n+1}=\tau_n\beta_i^n-\bsfp_i^{n+1}$.
\end{proof}
The viscous part of \eqref{eq:split} is already a sum of antisymmetric edge
fluxes, $\tau_n\big(d_{ij}^{{\rm L},n}-d_{ij}^{{\rm H},n}\big)\big(\bsfU_i^n-\bsfU_j^n\big)$,
so only the element vectors $\bsfa^{K,n+1}$ have to be converted into fluxes.

\subsection*{Element Neumann problems.}
We want to write $\widehat\bsfA_i^{n+1}$ as a flux over the edges. The viscous term in \eqref{eq:split} is already a sum of edge fluxes; Lemma~\ref{lem:zero_sum_K} shows that each $\bsfa_i^{K, n+1}$ has zero sum over $K$ for every $K\in\calT_h$; therefore, we can construct a local graph-Poisson problem to represent $\bsfa_i^{K, n+1}$ as a local antisymmetric flux. 
  
Let $L_K$ be the graph Laplacian of $G_K = \big(\calI(K), \calE(K) \big)$, $(L_K \bsfq^{K,n+1})_i=\sum_{j:\{i,j\}\in\calE(K)}(\bsfq_i^{K,n+1}-\bsfq_j^{K,n+1})$, applied componentwise. Since $G_K$ is connected, $\ker L_K$ consists of the constants, and by Lemma~\ref{lem:zero_sum_K}, the problem
\begin{equation}\label{eq:LK}
  L_K\bsfq^{K,n+1}=\bsfa^{K,n+1},\qquad\sum_{i\in\calI(K)}\bsfq_i^{K,n+1}=\boldsymbol 0,
\end{equation}
has a unique solution
\begin{equation}\label{eq:GK}
  \bsfq^{K,n+1}=G_k\,\bsfa^{K,n+1},\qquad G_k:=\big(L_K+\tfrac{1}{n_k}\mathbf 1\mathbf 1^\top\big)^{-1},
\end{equation}
where $n_k$ is the number of nodes in $K$. Indeed, $\mathbf 1^\top(L_K+\frac1{n_k}\mathbf 1\mathbf 1^\top)\bsfq^{K,n+1}=\mathbf 1^\top\bsfa^{K,n+1}=0$ gives $\mathbf 1^\top\bsfq^{K,n+1}=0$ and then $L_K\bsfq^{K,n+1}=\bsfa^{K,n+1}$. In addition, the matrix $G_k$ depends only on the polynomial degree and the numbering of the reference nodes; therefore, it is computed only once. By the minimal-norm property of the graph-Poisson fluxes \citep[Lemma~5.3]{Nazarov_2026}, $\bsfq_i^{K,n+1}-\bsfq_j^{K,n+1}$, $\{i,j\}\in\calE(K)$, are the fluxes of the smallest $\ell^2$-norm on $\calE(K)$ whose divergence is $\bsfa^{K,n+1}$. The same construction is used on a single element in
\citep[Lemma~4.2]{Abgrall_2018} to write residual distribution schemes in
flux form.

After solving \eqref{eq:LK} for every $K\in\calT_h$, we can finally write the edge fluxes as
\begin{equation}\label{eq:Aij_local}
  \widehat{\bsfA}_{ij}^{n+1}
  :=
  \sum_{K:\,\{i,j\}\in\calE(K)}\big(\bsfq_i^{K,n+1}-\bsfq_j^{K,n+1}\big)
  +\tau_n\big(d_{ij}^{{\rm L},n}-d_{ij}^{{\rm H},n}\big)\big(\bsfU_i^n-\bsfU_j^n\big),
  \qquad\{i,j\}\in\calE .
\end{equation}

Since $\sum_{j:\{i,j\}\in\calE(K)}\big(\bsfq_i^{K,n+1}-\bsfq_j^{K,n+1}\big)
=\big(L_K\bsfq^{K,n+1}\big)_i=\bsfa_i^{K,n+1}$, the fluxes \eqref{eq:Aij_local} are
antisymmetric and, by Lemma~\ref{lem:split},
\[
  \sum_{j\in\calN(i)}\widehat\bsfA_{ij}^{n+1}=\widehat\bsfA_i^{n+1},
  \qquad i\in\Vf .
\]

\section[Convex limiting]{Convex limiting}
\label{Sec:convex_limiting}
We combine the two forms \eqref{eq:HO_LO_pk} and \eqref{eq:HO_LO_p1} of Section~\ref{Sec:HO_LO_relation} into one formula,
\begin{equation}\label{eq:HO_LO}
  \UH_i = \ULtilde_i  +  \frac{1}{\mu_i} \sum_{j\in\calN(i)}\bsfA_{ij}^{n+1} + \iota \, \bsfs_i^{n+1},
\end{equation}
with $\mu_i=\mk_i$, $\iota=0$, $\ULtilde_i = \ULhat_i$ in the $\polP_k$-form (for $k=1,3$) and $\mu_i=\mone_i$,
$\iota=1$, $\ULtilde_i = \UL_i$ in the fine $\polP_1$-form (for $k=2,4$).

Our aim is to construct limiting coefficients such that the limited relation \eqref{eq:HO_LO} gives an admissible high-order solution.

\subsection{Local bounds}
For every node $i\in\calV$ we define bounds from the states that the low-order method combines, namely the neighbouring states and the bar states of \eqref{eq:LO_bar}:
\begin{align}
  \rho_i^{\min}
  &:=
  \min_{j\in\calN(i)}
  \big\{\rho_i^n,\ \overline{\rho}_{ij},\ \rho_j^n\big\},
  &
  \rho_i^{\max}
  &:=
  \max_{j\in\calN(i)}
  \big\{\rho_i^n,\ \overline{\rho}_{ij},\ \rho_j^n\big\},
  \label{eq:rho_minmax}
  \\
  s_i^{\min}
  &:=
  \min_{j\in\calN(i)}
  \big\{s(\bsfU_i^n),\ s(\Ubar_{ij}),\ s(\bsfU_j^n)\big\},
  &
  p_{\min}&:=10^{-14}.
  \label{eq:s_min}
\end{align}

The bar states are evaluated as in \eqref{eq:LO_bar}. The specific entropy of
a bar state is needed only when it can be smaller than 
$
\min_{j\in\calN(i)}
\big\{s(\bsfU_i^n),\ s(\bsfU_j^n)\big\}.
$
In addition we require
$p\ge p_{\min}$, that is,
$\rho E-\frac12|\bm|^2-\frac{p_{\min}}{\gamma-1}\rho\ge0$. This is the
internal-energy constraint of \citep{Guermond_Nazarov_Popov_Tomas_2018} with
the local bound replaced by a constant floor.

The base state of the fine $\polP_1$-form satisfies the above bounds by Proposition~\ref{prop:LO}; in the $\polP_k$-form, the low-order update satisfies Proposition~\ref{prop:LO} with a smaller time-step, see the discussion after \eqref{eq:HO_LO_pk}. To be safe against round-off and the wave-speed estimate, the bounds are enlarged to contain the base: $\rho_i^{\min}\leftarrow\min\{\rho_i^{\min},\widetilde{\rho\ }^{{\rm L}, n+1}_i\}$,
$\rho_i^{\max}\leftarrow\max\{\rho_i^{\max},\widetilde{\rho\ }^{{\rm L}, n+1}_i\}$,
$s_i^{\min}\leftarrow\min\{s_i^{\min},s(\ULtilde_i)-10^{-12}\}$.

It is well known that the strict local bounds clip smooth extrema, see \eg \citep{Guermond_Nazarov_Popov_Tomas_2018,Nazarov_2026} and the references therein. With the dimensionless
mesh size $\widetilde h_i:=(\mone_i/|\Omega|)^{1/d}$,
$r_i:=\min\{\frac12,C_{\rm rel}\widetilde h_i^{3/2}\}$,
$\delta_i:=C_{\rm rel}\widetilde h_i^{3/2}\rho_{\rm scale}$, where $\rho_{\rm scale}$ is
the maximum of the mean density and its oscillation over the mesh, we use
\begin{equation}\label{eq:relax_rho}
  \widecheck{\rho\ }_i^{\min}:=\max\big\{(1-r_i)\rho_i^{\min},\ \rho_i^{\min}-\delta_i\big\},
  \qquad
  \widecheck{\rho\ }_i^{\max}:=\rho_i^{\max}+\min\{r_i\rho_i^{\max},\ \delta_i\},
\end{equation}
and, with $r_i^s:=\min\{\frac12,C_{\rm rel,s} \widetilde h_i^{3/2}\}$,
\begin{equation}\label{eq:relax_s}
  \widecheck{s\ }_i^{\min}:=s_i^{\min}+\log\big(1-r_i^s\big).
\end{equation}

Every limited state then satisfies
$\rho\ge\widecheck{\rho\ }_i^{\min}\ge(1-r_i)\rho_i^{\min}\ge\frac12\rho_i^{\min}$
and $p=e^{s}\rho^\gamma\ge(1-r_i^s)e^{s_i^{\min}}\rho^\gamma
\ge\frac12e^{s_i^{\min}}\big(\frac12\rho_i^{\min}\big)^\gamma>0$; here
$\rho_i^{\min}>0$ and $s_i^{\min}>-\infty$ since $\bsfU^n$ and the bar states
are admissible. The factor $(1-r_i)$ is the positivity threshold of
\citep[Sec.~4.7]{Guermond_Nazarov_Popov_Tomas_2018}; unlike there, we always
relax the entropy bound by the full factor $(1-r_i^s)$. Over $N$ stages the
minimum of $e^{s}$ can therefore decrease at most by the factor
$(1-r^s)^N\approx e^{-Nr^s}$, which tends to one as $\widetilde h\to0$ on
quasi-uniform meshes, since $N\sim \widetilde h^{-1}$ and $r^s\sim \widetilde h^{3/2}$.

\subsection{Limited update}
\label{Sec:limited_update}
Let $\widehat\bsfA_{ij}^{n+1}$ be the fluxes of Section~\ref{Sec:global} or
Section~\ref{Sec:element_local}; at a free node they satisfy
$\sum_{j\in\calN(i)}\widehat\bsfA_{ij}^{n+1}=\widehat\bsfA_i^{n+1}$. Only the
fluxes on free edges are limited as edge fluxes: a prescribed node is
overwritten by the data and cannot take the opposite flux. The rest is
collected in the nodal source
\begin{equation}\label{eq:source}
  \bsfr_i^{n+1}:=\widehat\bsfA_i^{n+1}-\sum_{j\in\Nf(i)}\widehat\bsfA_{ij}^{n+1}-\bsfB_i^{n+1},
  \qquad i\in\Vf,
\end{equation}
that is, the flux $\sum_{b\in\calN(i)\cap\Vd}\widehat\bsfA_{ib}^{n+1}$ to the
prescribed neighbours minus the boundary term $\bsfB_i^{n+1}$ of
\eqref{eq:A_compatible}; at slip-wall nodes the density and energy components
of $-\bsfB_i^{n+1}$ are left out. At a far-field or open-boundary node $b$
of Section~\ref{Sec:boundary}, whose data are computed from $\UL_b$, the flux
$\widehat\bsfA_{ib}^{n+1}$ is dropped instead of being added to
$\bsfr_i^{n+1}$, so that the edge $\{i,b\}$ remains low order. Since
$\widehat\bsfA_i^{n+1}-\bsfB_i^{n+1}=\bsfA_i^{n+1}$, \eqref{eq:HO_LO} reads
\begin{equation}\label{eq:HO_LO_free}
  \UH_i=\ULtilde_i+\frac{1}{\mu_i}\Big(\sum_{j\in\Nf(i)}\widehat\bsfA_{ij}^{n+1}
  +\bsfr_i^{n+1}\Big)+\iota\,\bsfs_i^{n+1},\qquad i\in\Vf .
\end{equation}

We follow the limiting technique of \citep{Guermond_Nazarov_Popov_Tomas_2018}
and construct a limiting matrix $\big(\ell_{ij}\big)_{\{i,j\}\in\Ef}$ and
limiting vectors $(\eta_i)$, $(\zeta_i)$ with $0\le\ell_{ij}=\ell_{ji}\le1$ and
$0\le\eta_i,\zeta_i\le1$, so that the state
\begin{equation}\label{eq:lim}
  \V_i=\ULtilde_i+\frac{1}{\mu_i}\Big(\sum_{j\in\Nf(i)}\ell_{ij}\widehat\bsfA_{ij}^{n+1}
  +\zeta_i\bsfr_i^{n+1}\Big)+\iota\,\eta_i\bsfs_i^{n+1},
  \qquad i\in\Vf,
\end{equation}
stays admissible; at the prescribed nodes $\V_i$ is the boundary data. If all
limiting parameters vanish, $\V=\ULtilde$, the low-order base state; if they
are all one, $\V_i=\UH_i$ at every free node by \eqref{eq:HO_LO_free}, up to
the left-out wall components of $\bsfB_i^{n+1}$ and the dropped fluxes.

Let $n_i:=|\Nf(i)|+\iota+\kappa_i$ be the number of contributions at $i$,
where $\kappa_i=1$ if $\bsfr_i^{n+1}\neq\boldsymbol 0$ and $\kappa_i=0$
otherwise, and define the directions
\[
  \bsfP_{ij}^{n+1}:=\frac{n_i}{\mu_i}\widehat\bsfA_{ij}^{n+1},
  \qquad
  \bsfP_i^{\bsfr,n+1}:=\frac{n_i}{\mu_i}\bsfr_i^{n+1},
  \qquad
  \bsfP_i^{\bsfs,n+1}:=n_i\,\bsfs_i^{n+1}.
\]
Then \eqref{eq:lim} is the convex combination
\begin{equation}\label{eq:convcom}
  \V_i=\frac{1}{n_i}\Big(
  \sum_{j\in\Nf(i)}\big(\ULtilde_i+\ell_{ij}\bsfP_{ij}^{n+1}\big)
  +\kappa_i\big(\ULtilde_i+\zeta_i\bsfP_i^{\bsfr,n+1}\big)
  +\iota\big(\ULtilde_i+\eta_i\bsfP_i^{\bsfs,n+1}\big)\Big),
\end{equation}
and it suffices that every state in the parentheses lies in the convex set
\begin{equation}\label{eq:Bi}
  \calB_i:=\Big\{\bsfU\in\calA:\
  \widecheck{\rho\ }_i^{\min}\le\rho(\bsfU)\le\widecheck{\rho\ }_i^{\max},\
  p(\bsfU)\ge p_{\min},\ s(\bsfU)\ge \widecheck{s\ }_i^{\min}\Big\}.
\end{equation}

\subsection{Line search}
\label{Sec:line_search}
For a state $\bsfW_i\in\calB_i$ and a direction $\bsfP=(P^\rho,\bP^{\bm},P^E)$,
let $\Theta_i(\bsfW_i,\bsfP)$ be the largest $\theta\in[0,1]$ such that
$\bsfW_i+\chi \bsfP\in\calB_i$ for all $\chi\in[0,\theta]$, and write
$\rho(\theta)$, $\bm(\theta)$, $E(\theta)$ for the components of
$\bsfW_i+\theta\bsfP$. In the limiter, $\bsfW_i$ is the base state $\ULtilde_i$,
or the limited state of the previous pass.

Below we treat the three constraints in turn.

\paragraph{Density.} The density is affine in $\theta$, which gives
$\theta_\rho=\min\{1,(\widecheck{\rho\ }_i^{\max}-\rho_0)/P^\rho\}$ if $P^\rho>0$,
$\theta_\rho=\min\{1,(\widecheck{\rho\ }_i^{\min}-\rho_0)/P^\rho\}$ if $P^\rho<0$,
and $\theta_\rho=1$ otherwise; here $\rho_0:=\rho(0)$.

\paragraph{Pressure.} With $\e_{\min}:= p_{\min}/(\gamma - 1)$, the condition $p\ge p_{\min}$ for $\rho > 0$ is equivalent to
\begin{equation}
  \label{eq:pressure_g_theta}
  g(\theta) := \rho(\theta) \big( E(\theta) - \e_{\min} \big) - \frac12 |\bm(\theta)|^2 \ge 0,
\end{equation}
which is a quadratic polynomial with $g(0)>0$. If $g(\theta_\rho)\ge0$, we set $\theta_p:=\theta_\rho$; otherwise, $\theta_p$ is the smallest root of $g$ in $(0,\theta_\rho)$. Since $\{g\ge0\}$ is convex along the
segment, $g\ge0$ on $[0,\theta_p]$ \citep[Lemma~4.11]{Guermond_Nazarov_Popov_Tomas_2018}.

\paragraph{Specific entropy.} On $[0,\theta_p]$ we have $\rho(\theta)>0$ and
$p(\theta)>0$, and we use
\begin{equation}\label{eq:psi}
  \psi(\theta):=\rho(\theta)\big(s(\bsfW_i+\theta\bsfP)-\widecheck{s\ }_i^{\min}\big),
  \qquad\theta\in[0,\theta_p],
\end{equation}
which has the sign of $s-\widecheck{s\ }_i^{\min}$. It is concave, since $\rho s$
is a concave function of the conserved variables
\citep[Thm.~2.1]{Harten_Lax_Levermore_Morojoff_1998} and $\rho(\theta)$ is
affine. The base state is admissible, so $\psi(0)\ge0$. If
$\psi(\theta_p)\ge0$, we set $\theta_s:=\theta_p$. Otherwise, by concavity, the
set where $\psi\ge0$ is an interval $[0,\theta^*]$ with $\theta^*<\theta_p$.
We keep a bracket $[\theta^-,\theta^+]$ with $\psi(\theta^-)\ge0>\psi(\theta^+)$,
starting from $[0,\theta_p]$, and shrink it by Newton steps from $\theta^+$ and
secant steps from $\theta^-$, with bisection if a step leaves the bracket. We
stop when $\theta^+-\theta^-\le10^{-6}\theta_p$ and set $\theta_s:=\theta^-$,
which is admissible. Finally,
$\Theta_i(\bsfW_i,\bsfP):=(1-10^{-12})\,\theta_s$.

\begin{remark}[Avoiding powers in the line search]
  Every evaluation of $\psi$ needs $\log p$ and $\log\rho$. Since
  $s\ge\widecheck{s\ }_i^{\min}$ is equivalent to $p\ge\sigma_i\rho^\gamma$ with
  $\sigma_i:=e^{\widecheck{s\ }_i^{\min}}$, and $p=(\gamma-1)q/\rho$ with
  $q:=\rho E-\frac12|\bm|^2$, the condition reads
  $(\gamma-1)q\ge\sigma_i\rho^{\gamma+1}$. Along the segment $q(\theta)$ is a
  quadratic polynomial, and with $z(\theta):=\theta P^\rho/\rho(0)$ we have
  $\rho(\theta)=\rho(0)(1+z)$, so that the condition is
  \[
    (\gamma-1)\,q(\theta)\ge\sigma_i\rho(0)^\gamma\,\rho(\theta)\,(1+z)^\gamma .
  \]
  For $1\le\gamma\le2$ and $z>-1$, Taylor's formula with the Lagrange
  remainder gives
  \[
    (1+z)^\gamma\le{\rm UB}(z):=1+\gamma z+c_2z^2+c_3z^3+c_4z^4\max\{1,(1+z)^{-3}\},
    \qquad c_l=\binom{\gamma}{l},
  \]
  since $c_4\ge0$ and $(1+\xi)^{\gamma-4}\le\max\{1,(1+z)^{-3}\}$ for every
  $\xi$ between $0$ and $z$. Hence
  \[
    \psi_1(\theta):=(\gamma-1)\,q(\theta)-\sigma_i\rho(0)^\gamma\,\rho(\theta)\,
    {\rm UB}\big(z(\theta)\big)\ge0\quad\Longrightarrow\quad\psi(\theta)\ge0,
  \]
  so every $\theta$ accepted with $\psi_1$ is admissible. $\psi_1$ needs no
  logarithm, and $\rho(0)^\gamma$ is computed once per node and pass for all
  directions. We use $\psi_1$ when $|z(\theta_p)|\le0.1$, where for $\gamma=1.4$
  ${\rm UB}$ exceeds $(1+z)^\gamma$ by less than $10^{-6}$ relative, and $\psi$
  otherwise.
\end{remark}

\paragraph{Limiters.}
Once the line search is done, we define the symmetric edge limiters and the source limiters as
\begin{equation}\label{eq:euler_edge_limiter}
  \begin{cases}
    \ell_{ij}:=\min\big\{\Theta_i(\ULtilde_i,\bsfP_{ij}^{n+1}),\ \Theta_j(\ULtilde_j,\bsfP_{ji}^{n+1})\big\},
    \\
    \eta_i:=\Theta_i(\ULtilde_i,\bsfP_i^{\bsfs,n+1}),
    \\
  \zeta_i:=\Theta_i(\ULtilde_i,\bsfP_i^{\bsfr,n+1}).
\end{cases}
\end{equation}

We follow \citep{Nazarov_2026} and apply the limiter in passes. Let
$\bsfV^{[0]}:=\ULtilde$, $\widehat\bsfA_{ij}^{[1]}:=\widehat\bsfA_{ij}^{n+1}$,
$\bsfs_i^{[1]}:=\bsfs_i^{n+1}$ and $\bsfr_i^{[1]}:=\bsfr_i^{n+1}$. Pass $\nu\ge1$
computes the limiters \eqref{eq:euler_edge_limiter} for these parts, with
$\bsfV^{[\nu-1]}$ in place of $\ULtilde$ and with the same bounds and weights,
sets
\[
  \bsfV_i^{[\nu]}:=\bsfV_i^{[\nu-1]}+\frac{1}{\mu_i}\Big(\sum_{j\in\Nf(i)}\ell_{ij}^{[\nu]}\widehat\bsfA_{ij}^{[\nu]}
  +\zeta_i^{[\nu]}\bsfr_i^{[\nu]}\Big)+\iota\,\eta_i^{[\nu]}\bsfs_i^{[\nu]},
\]
and passes on the rejected parts
$\widehat\bsfA_{ij}^{[\nu+1]}:=(1-\ell_{ij}^{[\nu]})\widehat\bsfA_{ij}^{[\nu]}$,
$\bsfs_i^{[\nu+1]}:=(1-\eta_i^{[\nu]})\bsfs_i^{[\nu]}$ and
$\bsfr_i^{[\nu+1]}:=(1-\zeta_i^{[\nu]})\bsfr_i^{[\nu]}$. After $\nu_{\max}$
passes, $\V:=\bsfV^{[\nu_{\max}]}$ is \eqref{eq:lim} with the accepted fractions
\[
  \bar\ell_{ij}:=1-\prod_{\nu=1}^{\nu_{\max}}\big(1-\ell_{ij}^{[\nu]}\big),\qquad
  \bar\eta_i:=1-\prod_{\nu=1}^{\nu_{\max}}\big(1-\eta_i^{[\nu]}\big),\qquad
  \bar\zeta_i:=1-\prod_{\nu=1}^{\nu_{\max}}\big(1-\zeta_i^{[\nu]}\big)
\]
in place of $\ell_{ij}$, $\eta_i$, $\zeta_i$; they lie in $[0,1]$, and
$\bar\ell_{ij}=\bar\ell_{ji}$. We use $\nu_{\max}=3$ and stop earlier when
nothing was limited in a pass.

\subsection{Properties of the limited solution}

Every pass is a limited update of the form \eqref{eq:lim}, with
$\bsfV^{[\nu-1]}$ in place of $\ULtilde$ and the rejected parts in place of
$\widehat\bsfA^{n+1}$, $\bsfs^{n+1}$, $\bsfr^{n+1}$, and after $\nu_{\max}$
passes $\V$ is \eqref{eq:lim} with $\bar\ell_{ij}$, $\bar\eta_i$, $\bar\zeta_i$.
It therefore suffices to state the properties for one pass:
Theorem~\ref{thm:euler_IDP} applies to every pass, with base state
$\bsfV^{[\nu-1]}\in\calB_i$, and Theorem~\ref{thm:conservation} and the
exchange of Section~\ref{Sec:exchange}, which is applied once after the last
pass, hold with $\bar\ell_{ij}$, $\bar\eta_i$, $\bar\zeta_i$ in place of
$\ell_{ij}$, $\eta_i$, $\zeta_i$.

\begin{theorem}[Invariant domain]\label{thm:euler_IDP}
  Assume that $\ULtilde_i\in\calB_i$ for all $i\in\Vf$ and that the boundary
  data lie in $\calA$. Then the limited state \eqref{eq:lim} with the limiters
  \eqref{eq:euler_edge_limiter} satisfies $\V_i\in\calB_i\subset\calA$ for all
  $i\in\Vf$ and $\V_i\in\calA$ for all $i\in\Vd$. The same holds after the
  slip projection, and every step of the strong stability preserving
  Runge--Kutta method stays in $\calA$.
\end{theorem}
\begin{proof}
  The set $\calB_i$ is convex: the density bounds define a slab, $p\ge p_{\min}$
  is $\rho\varepsilon(\bsfU)\ge p_{\min}/(\gamma-1)$ with the concave function $\rho\varepsilon$,
  and $\{s\ge\widecheck{s\ }_i^{\min}\}$ is convex since $s$ is quasiconcave.
  By the definition of $\Theta_i$, the segment from $\ULtilde_i$ to
  $\ULtilde_i+\Theta_i(\ULtilde_i,\bsfP)\bsfP$ lies in $\calB_i$. Since
  $\ell_{ij}\le\Theta_i(\ULtilde_i,\bsfP_{ij}^{n+1})$ by
  \eqref{eq:euler_edge_limiter}, the state $\ULtilde_i+\ell_{ij}\bsfP_{ij}^{n+1}$
  lies on this segment and hence in $\calB_i$; so do
  $\ULtilde_i+\eta_i\bsfP_i^{\bsfs,n+1}$ and $\ULtilde_i+\zeta_i\bsfP_i^{\bsfr,n+1}$.
  By \eqref{eq:convcom}, $\V_i$ is a convex combination of these states, so
  $\V_i\in\calB_i$. At a prescribed node $\V_i$ is the boundary data. The slip
  projection keeps $\rho_i$ and $E_i$ and does not increase $|\bm_i|$, so it
  does not decrease $p$ or $s$ (Section~\ref{Sec:boundary}). Finally, $\calA$
  is convex and a step of the strong stability preserving Runge--Kutta method
  is a convex combination of forward Euler stages.
\end{proof}

\begin{theorem}[Conservation]\label{thm:conservation}
  The limited state \eqref{eq:lim} satisfies
  \begin{equation}\label{eq:cons_node}
    \mu_i\big(\V_i-\bsfU_i^n\big)=\mone_i\big(\UL_i-\bsfU_i^n\big)
    +\sum_{j\in\Nf(i)}\ell_{ij}\widehat\bsfA_{ij}^{n+1}
    +\zeta_i\bsfr_i^{n+1}+\iota\,\eta_i\,\mone_i\bsfs_i^{n+1},
    \qquad i\in\Vf,
  \end{equation}
  and hence
  \begin{equation}\label{eq:cons_sum}
    \sum_{i\in\Vf}\mu_i\big(\V_i-\bsfU_i^n\big)
    =\sum_{i\in\Vf}\mone_i\big(\UL_i-\bsfU_i^n\big)
    +\sum_{i\in\Vf}\big(\zeta_i\bsfr_i^{n+1}+\iota\,\eta_i\,\mone_i\bsfs_i^{n+1}\big).
  \end{equation}
  For periodic boundaries and exact linear solves, $\Mk(\V)=\Mk(\bsfU^n)$ in
  the $\polP_k$-form, and in the fine $\polP_1$-form
  \begin{equation}\label{eq:mass_defect}
    \Mk(\V)-\Mk(\bsfU^n)=\sum_{K\in\calT_h}\sum_{i\in\calI(K)}
    \big(m_{i,K}^{\polP_1}-m_{i,K}^{\polP_k}\big)\bsfW_i,
  \end{equation}
  where
  \begin{equation}\label{eq:Wi}
    \bsfW_i:=\eta_i\big(\UH_i-\bsfU_i^n\big)-\big(\V_i-\bsfU_i^n\big),
  \end{equation}
    with $\bsfW_i=\boldsymbol 0$ at every node where nothing is limited.
\end{theorem}
\begin{proof}
  Multiplying \eqref{eq:lim} by $\mu_i$ and using
  $\mu_i(\ULtilde_i-\bsfU_i^n)=\mone_i(\UL_i-\bsfU_i^n)$ in both forms, and
  $\mu_i=\mone_i$ when $\iota=1$, gives \eqref{eq:cons_node}. In the sum over
  $\Vf$ the edge terms cancel, since $\ell_{ij}=\ell_{ji}$ and
  $\widehat\bsfA_{ji}^{n+1}=-\widehat\bsfA_{ij}^{n+1}$.
  With periodic boundaries
  all nodes are free and $\bsfB^{n+1}=\boldsymbol 0$, and if the mass system \eqref{eq:HO} and, for the
  global fluxes, the graph-Poisson system \eqref{eq:L_dirichlet} are solved
  exactly
  $\bsfr^{n+1}=\boldsymbol 0$ by \eqref{eq:source}; moreover
  $\Mone(\UL)=\Mone(\bsfU^n)$ by Proposition~\ref{prop:LO}. In the $\polP_k$-form
  $\mu=\mk$ and $\iota=0$, which gives $\Mk(\V)=\Mk(\bsfU^n)$. In the fine
  $\polP_1$-form, $\mone_i\bsfs_i^{n+1}=(\mone_i-\mk_i)(\UH_i-\bsfU_i^n)$, so
  $\Mone(\V)-\Mone(\bsfU^n)=\sum_i\eta_i(\mone_i-\mk_i)(\UH_i-\bsfU_i^n)$;
  subtracting $\sum_i(\mone_i-\mk_i)(\V_i-\bsfU_i^n)$ and splitting
  $\mone_i-\mk_i$ over the elements gives the defect. If nothing is limited
  at $i$, then $\V_i=\UH_i$ and $\eta_i=1$, so $\bsfW_i=\boldsymbol 0$.
\end{proof}

The defect is thus confined to the elements $K$ with a limited node, \ie
with $\ell_{ij}<1$ for some edge at a node of $K$, or $\eta_i<1$ or
$\zeta_i<1$ at a node of $K$. It is removed element by element, without
violating the bounds, by the exchange of Section~\ref{Sec:exchange}.

\subsection{Element-local exchange for $\polP_2$ and $\polP_4$}
\label{Sec:exchange}

Redistributing a conservation defect among nearby states within their
bounds is the idea of the ``repair'' methods
\citep{Kucharik_Shashkov_Wendroff_2003,Shashkov_Wendroff_2004,Loubere_Staley_Wendroff_2006}
and of conservative redistribution limiting
\citep{Guermond_Wang_2025,Guermond_2026}; an element-wise mass correction is
used in \citep{Anderson_2017}, and distributing a defect in proportion to the
available room is the clip-and-assured-sum fixer of \citep{Bradley_2019}.
Here the states are already admissible and the defect stems from the two
masses; the line search measures the room in the convex set $\calB_i$, and,
as in the cell-based algorithm of \citep{Guermond_2026}, all elements use the
same states $\V$, so that the result does not depend on the order of the
elements.

In the fine $\polP_1$-form the limited state conserves $\Mone$ up to the
limited sources, but not $\Mk$. By \eqref{eq:mass_defect}, the $\Mk$ defect is
the sum of the element defects
\[
  D_K:=\sum_{i\in\calI(K){\cap\Vf}}\big(m_{i,K}^{\polP_1}-m_{i,K}^{\polP_k}\big)\bsfW_i,
\]
with $\bsfW_i$ of \eqref{eq:Wi}, where $m_{i,K}^{\polP_1}$ and
$m_{i,K}^{\polP_k}$ are the contributions of $K$ to $\mone_i$ and $\mk_i$;
$D_K$ vanishes on every element where nothing is limited.
The prescribed
nodes are left out: their state is overwritten by the boundary data, so
$\bsfW_i$ is not a defect of the update there; as in the $\polP_k$-form, $\UH_i$
is not used at these nodes. With periodic boundaries all nodes are free, and
$\sum_{K\in \calT_h}D_K$ is the defect \eqref{eq:mass_defect}. As in
\citep{Nazarov_2026}, we remove it element by element by increments at the
free nodes of $K$ with $\mk_i>0$. Since the admissible set of the Euler
system is not a box, the capacities of \citep{Nazarov_2026} are replaced by
the line search of Section~\ref{Sec:line_search}. Let
$\calT(i):=\{K\in\calT_h:\ i\in\calI(K)\}$. For every element $K$:
\begin{enumerate}
\item set $\bsft_K:=-D_K$;
\item for every free node $i\in\calI(K)$ with $\mk_i>0$, compute
  $\vartheta_i^K:=\Theta_i(\V_i,\bsfP_i^K)$ with the direction
  $\bsfP_i^K:=|\calT(i)|\,\bsft_K/\mk_i$, and set
  $S_K:=\sum_{i\in\calI(K)}\vartheta_i^K$ and $\alpha_K:=\min\{1,S_K\}$;
\item if $S_K>0$, add
  $\delta\bsfU_i^K:=\frac{\alpha_K\vartheta_i^K}{S_K}\,\frac{\bsft_K}{\mk_i}$
  to every such node.
\end{enumerate}
The increments of $K$ change $\Mk$ by
$\sum_{i\in\calI(K)}\mk_i\,\delta\bsfU_i^K=\alpha_K\bsft_K$, so they remove the whole defect
of $K$ if $\alpha_K=1$. Since
$\delta\bsfU_i^K=\frac{\alpha_K}{S_K}\frac{1}{|\calT(i)|}\,\vartheta_i^K\bsfP_i^K$
with $\alpha_K\le S_K$, the final state $\V_i+\sum_{K\in\calT(i)}\delta\bsfU_i^K$
is a convex combination of $\V_i$ and the states
$\V_i+\vartheta_i^K\bsfP_i^K\in\calB_i$, so the exchange preserves the bounds.
If $\alpha_K<1$, the nodes of $K$ have no room left in the direction of
$\bsft_K$, and the part $(1-\alpha_K)\bsft_K$ is not restored.

\paragraph{Second round.} This happens on few elements, typically next to strong shocks,
and their neighbours usually still have room. We therefore offer the leftover
$\bsft'_K:=(1-\alpha_K)\bsft_K$ to the patch
$\calP_K:=\bigcup\{\calI(K'):\ K'=K\text{ or }K'\text{ shares a facet with }K\}$,
by the same three steps with three changes: the base state is
${\V_i}':=\V_i+\sum_{K\in\calT(i)}\delta\bsfU_i^K$, the receiving nodes are the
free nodes of $\calP_K$ with $\mk_i>0$, and $|\calT(i)|$ is replaced by the
number of patches that contain $i$. The convexity argument above applies
with ${\V_i}'$ in place of $\V_i$, so the bounds still hold.

\paragraph{Further rounds.}
In the rare cells without room, the patch can be grown for a few more rounds, followed by one global round if needed, as in \citep[\S6.1]{Loubere_Staley_Wendroff_2006} and \citep[Sec.~4]{Guermond_2026}. Each round keeps the bounds by the convexity argument above; the global round restores $\Mk$ exactly whenever the free nodes have enough room in total, but it conserves only globally.

In our simulations we use two rounds. In the smooth tests, the element round already restores the whole defect. With shocks, the second round restores most of what the element round leaves, and the rest is a small fraction of the total mass; it is largest for $\polP_2$, whose vertices have $\mk_i=0$ and cannot receive.

\section{Numerical validation }
\label{Sec:numerics}

For the smooth tests, we use the fourth-order five-stage strong stability
preserving Runge--Kutta method SSP-RK(5,4) of \citep{Kraaijevanger_1991}. For the remaining
discontinuous problems, we use the third-order SSP-RK(3,3) method of
\citep{Shu_Osher1988}. Both methods are convex combinations of forward Euler steps.

The time step is computed at every time level by
\begin{equation}\label{eq:dt}
  \tau_n=n_{\rm st} \textrm{ CFL } \frac{1}{r_{\max}}\,\min_{i\in\calV}\frac{\mone_i}{2\sum_{j\in\calN(i)}d_{ij}^{\rm L}(\bsfU^n)},
\end{equation}
where $n_{\rm st}$ is the number of stages, $n_{\rm st}=5$ for SSP-RK(5,4) and $n_{\rm st}=3$ for SSP-RK(3,3), and $\textrm{CFL}>0$ is a CFL number,
$r_{\max}=\max_{i\in\calV}\mone_i/\mk_i$ in the $\polP_k$ form ($r_{\max}=1$ for
$\polP_1$ and in the fine $\polP_1$-form, $r_{\max}=1.48$ for $\polP_3$). 

Algorithm~\ref{alg:stage} summarizes one forward Euler stage.
The method is implemented in C++ using the finite element library DOLFINx \citep{BarattaEtal2023}.

In the convergence tables we report the combined relative error
\begin{equation}\label{eq:delta_q}
  \delta_q
  :=
  \frac{\|\rho_h-\rho\|_{L^q(\Omega)}}{\|\rho\|_{L^q(\Omega)}}
  +
  \frac{\|\bm_h-\bm\|_{L^q(\Omega)}}{\|\bm\|_{L^q(\Omega)}}
  +
  \frac{\|E_h-E\|_{L^q(\Omega)}}{\|E\|_{L^q(\Omega)}},
  \qquad
  q=1,2,\infty .
\end{equation}
In the implementation, the exact solution is interpolated to $\polP_{k+2}$, then
the $L^1$ and $L^2$ norms in \eqref{eq:delta_q} are computed by quadrature on the elements, exactly for $q=2$, and the $L^\infty$ norm is the maximum over the nodes of the $\polP_{k+2}$ space.

Unless stated otherwise, in all simulations we use
$\gamma = 1.4$,
$C_{\textrm{R}} = 0.5$,
$C_{\max} = 1$, $c_{\textrm{cip}}=0.05$,
$C_{\textrm{rel}} = C_{\textrm{rel},s} = 1$.
The CFL number for all simulations is $\textrm{CFL}=0.3$.
All meshes are unstructured Delaunay triangulations. 

\begin{algorithm}[htbp]
  \caption{One forward Euler stage with graph-Poisson limiting}
  \label{alg:stage}
  \begin{algorithmic}[1]
    \Require $\bsfU^n$ (boundary conditions applied), $\tau_n$.
    \State Low-order state $\UL$ by \eqref{eq:LO}; boundary conditions.
    \State High-order state $\UH$ by \eqref{eq:HO}; boundary conditions.
    \State Local bounds \eqref{eq:rho_minmax}--\eqref{eq:relax_s} and base state
      $\ULtilde$ of \eqref{eq:HO_LO}. If $\UH_i\in\calB_i$ for all $i\in\Vf$,
      \Return $\V:=\UH$.
    \State Fluxes $\widehat\bsfA_{ij}^{n+1}$, global (Section~\ref{Sec:global}) or
      element-local (Section~\ref{Sec:element_local}), and sources
      $\bsfr_i^{n+1}$, $\bsfs_i^{n+1}$.
    \State Limited state $\V$ by \eqref{eq:lim}, with $\nu_{\max}$ passes of
      the limiter \eqref{eq:euler_edge_limiter}.
    \State If $\iota=1$, element-local exchange (Section~\ref{Sec:exchange}).
    \State Boundary conditions on $\V$.
    \Ensure $\V$.
  \end{algorithmic}
\end{algorithm}

\subsection{Isentropic vortex}
\label{Sec:vortex}
To validate the accuracy for a smooth solution, we consider the isentropic vortex
\[
  \bv=\bv_\infty+\frac{\beta}{2\pi}e^{(1-r^2)/2}(-\bar y,\bar x),\qquad
  \rho=\Big(1-\frac{(\gamma-1)\beta^2}{8\gamma\pi^2}e^{1-r^2}\Big)^{\frac1{\gamma-1}},\qquad
  p=\rho^\gamma,
\]
with $(\bar x,\bar y)=\bx-\bx_0-\bv_\infty t$, $r^2=\bar x^2+\bar y^2$,
$\beta=2$, $\bv_\infty=(1,1)$, $\bx_0=(-1,-1)$, is advected in
$\Omega=(-10,10)^2$ until $\widehat t=2$, with the exact solution as Dirichlet data on
$\p\Omega$. We use SSP-RK(5,4) 
on a sequence of meshes with comparable numbers of nodes for all polynomial degrees. We present the results in Table~\ref{tab:vortex_unstructured}. We observe
near-optimal rates in all norms.  

\begin{table}[h!]
  \caption{Isentropic vortex, $\widehat t=2$, unstructured meshes: combined relative
    errors $\delta_q$ and observed rates with $h=(\# \textrm{nodes})^{-1/2}$.}
  \label{tab:vortex_unstructured}
  \centering
  \scriptsize
  \begin{tabular}{lrcccccc}
    \toprule
    & $\# \textrm{nodes}$ & $\delta_1$ & rate & $\delta_2$ & rate & $\delta_\infty$ & rate\\
    \midrule
    $\polP_1$ & 1\,576 & $1.98\times10^{-3}$ & -- & $1.06\times10^{-2}$ & -- & $1.41\times10^{-1}$ & --\\
     & 4\,998 & $4.56\times10^{-4}$ & 2.54 & $2.52\times10^{-3}$ & 2.48 & $3.69\times10^{-2}$ & 2.32\\
     & 16\,158 & $1.07\times10^{-4}$ & 2.48 & $5.87\times10^{-4}$ & 2.49 & $9.02\times10^{-3}$ & 2.40\\
     & 50\,464 & $2.79\times10^{-5}$ & 2.35 & $1.53\times10^{-4}$ & 2.36 & $2.89\times10^{-3}$ & 2.00\\
     & 159\,544 & $7.75\times10^{-6}$ & 2.23 & $4.18\times10^{-5}$ & 2.25 & $8.00\times10^{-4}$ & 2.23\\
     & 499\,990 & $2.31\times10^{-6}$ & 2.12 & $1.25\times10^{-5}$ & 2.12 & $2.85\times10^{-4}$ & 1.81\\
    \midrule
    $\polP_2$ & 1\,613 & $1.15\times10^{-3}$ & -- & $5.75\times10^{-3}$ & -- & $8.43\times10^{-2}$ & --\\
     & 4\,841 & $1.38\times10^{-4}$ & 3.86 & $6.72\times10^{-4}$ & 3.91 & $9.43\times10^{-3}$ & 3.99\\
     & 15\,597 & $1.99\times10^{-5}$ & 3.31 & $9.76\times10^{-5}$ & 3.30 & $2.07\times10^{-3}$ & 2.59\\
     & 50\,233 & $3.24\times10^{-6}$ & 3.10 & $1.57\times10^{-5}$ & 3.13 & $3.79\times10^{-4}$ & 2.90\\
     & 159\,997 & $5.90\times10^{-7}$ & 2.94 & $2.92\times10^{-6}$ & 2.90 & $1.07\times10^{-4}$ & 2.18\\
     & 501\,117 & $1.05\times10^{-7}$ & 3.02 & $5.28\times10^{-7}$ & 3.00 & $1.87\times10^{-5}$ & 3.06\\
    \midrule
    $\polP_3$ & 1\,381 & $6.96\times10^{-4}$ & -- & $3.11\times10^{-3}$ & -- & $4.74\times10^{-2}$ & --\\
     & 4\,765 & $6.77\times10^{-5}$ & 3.76 & $3.43\times10^{-4}$ & 3.56 & $6.42\times10^{-3}$ & 3.23\\
     & 16\,495 & $4.88\times10^{-6}$ & 4.24 & $2.46\times10^{-5}$ & 4.25 & $9.22\times10^{-4}$ & 3.12\\
     & 50\,485 & $4.76\times10^{-7}$ & 4.16 & $2.21\times10^{-6}$ & 4.31 & $8.97\times10^{-5}$ & 4.17\\
     & 159\,421 & $5.36\times10^{-8}$ & 3.80 & $2.54\times10^{-7}$ & 3.76 & $9.42\times10^{-6}$ & 3.92\\
     & 495\,874 & $6.02\times10^{-9}$ & 3.85 & $2.98\times10^{-8}$ & 3.77 & $1.52\times10^{-6}$ & 3.22\\
    \midrule
    $\polP_4$ & 1\,537 & $5.16\times10^{-4}$ & -- & $1.83\times10^{-3}$ & -- & $2.17\times10^{-2}$ & --\\
     & 4\,841 & $2.41\times10^{-5}$ & 5.34 & $1.00\times10^{-4}$ & 5.07 & $2.33\times10^{-3}$ & 3.89\\
     & 15\,609 & $1.85\times10^{-6}$ & 4.39 & $9.06\times10^{-6}$ & 4.11 & $2.84\times10^{-4}$ & 3.59\\
     & 49\,665 & $1.12\times10^{-7}$ & 4.85 & $5.35\times10^{-7}$ & 4.89 & $1.52\times10^{-5}$ & 5.07\\
     & 161\,481 & $6.93\times10^{-9}$ & 4.72 & $3.50\times10^{-8}$ & 4.63 & $1.51\times10^{-6}$ & 3.91\\
     & 503\,969 & $5.43\times10^{-10}$ & 4.48 & $2.76\times10^{-9}$ & 4.46 & $1.29\times10^{-7}$ & 4.32\\
    \bottomrule
  \end{tabular}
\end{table}

\subsection{Isentropic vortex leaving the domain}
\label{Sec:isentropic}
To test the treatment of the Dirichlet boundaries, the vortex
from Section~\ref{Sec:vortex}, with the speed
$\bv_\infty=(1,1)$, $\beta=2$, and $\bx_0=(1,1)$
is advected in $\Omega = (-4,4)^2$ until $\widehat t = 3$.
At this final time, the center of the vortex lies on the corner
$(4,4)$. The exact solution is prescribed at all boundary
nodes, so $\Vd = \calV_{\p}$. As in the previous test, we use
SSP-RK(5,4). 

With the fluxes to the prescribed nodes kept, the convex limiting
with element-local fluxes of Section~\ref{Sec:element_local}
gives near-optimal rates for all polynomial degrees, see Table~\ref{tab:vortex_exit}. 

\begin{table}[h!]
  \caption{Vortex leaving the domain, $\widehat t=3$, unstructured meshes, element-local fluxes: combined relative
    errors $\delta_q$ and observed rates.}
  \label{tab:vortex_exit}
  \centering
  \scriptsize
  \begin{tabular}{lrcccccc}
    \toprule
    & $\# \textrm{nodes}$ & $\delta_1$ & rate & $\delta_2$ & rate & $\delta_\infty$ & rate\\
    \midrule
    $\polP_1$ & 1\,664 & $3.66\times10^{-4}$ & -- & $1.55\times10^{-3}$ & -- & $1.87\times10^{-2}$ & --\\
     & 4\,995 & $9.08\times10^{-5}$ & 2.53 & $3.73\times10^{-4}$ & 2.59 & $4.73\times10^{-3}$ & 2.50\\
     & 15\,915 & $2.59\times10^{-5}$ & 2.17 & $1.00\times10^{-4}$ & 2.27 & $1.57\times10^{-3}$ & 1.91\\
     & 50\,039 & $7.50\times10^{-6}$ & 2.16 & $2.88\times10^{-5}$ & 2.17 & $4.70\times10^{-4}$ & 2.10\\
     & 159\,660 & $2.27\times10^{-6}$ & 2.06 & $8.54\times10^{-6}$ & 2.09 & $1.41\times10^{-4}$ & 2.08\\
     & 501\,275 & $7.11\times10^{-7}$ & 2.03 & $2.68\times10^{-6}$ & 2.02 & $5.57\times10^{-5}$ & 1.62\\
    \midrule
    $\polP_2$ & 1\,649 & $8.74\times10^{-5}$ & -- & $3.20\times10^{-4}$ & -- & $4.32\times10^{-3}$ & --\\
     & 5\,017 & $1.63\times10^{-5}$ & 3.02 & $6.21\times10^{-5}$ & 2.95 & $9.90\times10^{-4}$ & 2.65\\
     & 15\,713 & $3.09\times10^{-6}$ & 2.91 & $1.26\times10^{-5}$ & 2.79 & $3.06\times10^{-4}$ & 2.06\\
     & 50\,125 & $5.20\times10^{-7}$ & 3.07 & $2.11\times10^{-6}$ & 3.08 & $6.97\times10^{-5}$ & 2.55\\
     & 158\,765 & $9.82\times10^{-8}$ & 2.89 & $4.21\times10^{-7}$ & 2.80 & $1.70\times10^{-5}$ & 2.45\\
     & 502\,677 & $1.74\times10^{-8}$ & 3.00 & $7.42\times10^{-8}$ & 3.01 & $3.87\times10^{-6}$ & 2.57\\
    \midrule
    $\polP_3$ & 1\,765 & $3.37\times10^{-5}$ & -- & $1.11\times10^{-4}$ & -- & $1.75\times10^{-3}$ & --\\
     & 5\,212 & $3.27\times10^{-6}$ & 4.31 & $1.09\times10^{-5}$ & 4.28 & $2.08\times10^{-4}$ & 3.94\\
     & 15\,580 & $3.86\times10^{-7}$ & 3.91 & $1.21\times10^{-6}$ & 4.01 & $2.63\times10^{-5}$ & 3.77\\
     & 50\,116 & $3.98\times10^{-8}$ & 3.89 & $1.26\times10^{-7}$ & 3.88 & $2.80\times10^{-6}$ & 3.84\\
     & 160\,195 & $4.20\times10^{-9}$ & 3.87 & $1.34\times10^{-8}$ & 3.85 & $3.52\times10^{-7}$ & 3.57\\
     & 500\,767 & $4.83\times10^{-10}$ & 3.80 & $1.59\times10^{-9}$ & 3.74 & $4.33\times10^{-8}$ & 3.68\\
    \midrule
    $\polP_4$ & 1\,537 & $1.62\times10^{-5}$ & -- & $4.61\times10^{-5}$ & -- & $5.56\times10^{-4}$ & --\\
     & 5\,521 & $8.27\times10^{-7}$ & 4.66 & $2.35\times10^{-6}$ & 4.66 & $5.68\times10^{-5}$ & 3.57\\
     & 16\,721 & $7.14\times10^{-8}$ & 4.42 & $2.19\times10^{-7}$ & 4.28 & $3.95\times10^{-6}$ & 4.81\\
     & 49\,953 & $5.64\times10^{-9}$ & 4.64 & $1.65\times10^{-8}$ & 4.73 & $2.67\times10^{-7}$ & 4.92\\
     & 158\,849 & $4.46\times10^{-10}$ & 4.39 & $1.30\times10^{-9}$ & 4.38 & $2.44\times10^{-8}$ & 4.14\\
     & 497\,609 & $3.64\times10^{-11}$ & 4.39 & $1.01\times10^{-10}$ & 4.48 & $2.70\times10^{-9}$ & 3.86\\
    \bottomrule
  \end{tabular}
\end{table}

\paragraph{Global fluxes.}
With the global fluxes of Section~\ref{Sec:global}, the isentropic vortex
converges at the same rates as in Table~\ref{tab:vortex_unstructured}: with the
limiter active in every stage, the errors agree with those of the
element-local fluxes to four digits on all meshes.

For the vortex leaving the domain, the global and the element-local fluxes
give the same errors, to three digits, on the coarser meshes.
On the finer meshes, however, the global fluxes lose the convergence rate. The loss starts
when the limiter begins to cut the global fluxes, and the error is largest
next to the outflow boundary. The reason is the nonlocality of the global fluxes.
In the free stream the density takes its global maximum $\rho=1$ and the entropy
is constant, so the bounds leave only the relaxation \eqref{eq:relax_rho}--\eqref{eq:relax_s}
as room. In the isentropic vortex the defect stays far from the boundary, and the
global fluxes it produces in this region are small. A similar loss of
accuracy for flat solutions was observed for scalar conservation
laws in \citep{Nazarov_2026}.

\subsection{Sod shock tube} \label{Sec:sod}
Now, we solve a discontinuous problem, namely the Sod problem. The computational domain
$\Omega = (0,1)\times (0, 0.2)$ is triangulated with unstructured meshes. The exact
solution is imposed on the left and right boundaries, and the slip boundary condition
is imposed on the top and bottom boundaries. The initial condition consists of
$(\rho, \bv, p) = (1, {\bold 0}, 1)$ for $x<0.5$, and 
$(\rho, \bv, p) = (0.125, {\bold 0}, 0.1)$ for $x>0.5$. To avoid interpolating a
discontinuity into $\polP_k$, the initial jump is replaced by
$
\rho_h^0 = \rho_R + (\rho_L - \rho_R)\theta_h
$
and
$
p_h^0 = p_R + (p_L - p_R)\theta_h,
$
with
$\theta_h(x)=\frac12\big(1-\tanh\frac{x-0.5}{\varepsilon_h}\big)$ and
$\varepsilon_h=0.2h/k$, where $h$ is the mean over $K\in\calT_h$ of the longest edge of $K$, so that $h/k$ is a typical node spacing. Note that the smoothing vanishes as $h \to 0$
with first order rate.

The simulation is run using SSP-RK(3,3).
We report the relative $L^1$ error of the density for all polynomial
degrees in Table~\ref{tab:sod}. The rates are as expected for
problems with shocks and contact discontinuities. The errors are similar for all degrees and slightly larger for $\polP_2$ and $\polP_4$; we suspect that interpolating the exact solution into $\polP_{k+2}$ for the $L^1$ norm produces oscillations at the discontinuities, which increase the error.

We also observe that the global flux limiter cuts 43--71\% of the edges,
and this grows under refinement. Since the global flux is nonlocal, the defect
produced at the shock and contact discontinuity spreads over the whole domain,
including the constant states. Moreover, around the constant states
the bounds leave almost no room, so even tiny fluxes get cut nearly
everywhere. The element-local fluxes stay where the defect is
generated, \ie near the shock and the contact discontinuity. When the mesh
is refined, these regions get smaller, and the share of limited edges decreases.

\begin{table}[h!]
  \caption{Sod problem, $\widehat t=0.2$: relative $L^1$ errors of the density,
    observed rates, and the fraction of the edges
    whose flux is limited in the first pass, for the global and the
    element-local fluxes.}
  \label{tab:sod}
  \centering
  \scriptsize
  \begin{tabular}{rcccccc}
    \toprule
    & \multicolumn{3}{c}{global} & \multicolumn{3}{c}{local}\\
    \cmidrule(lr){2-4}\cmidrule(lr){5-7}

    \#nodes & $L^1$ & rate & limited & $L^1$ & rate & limited\\
    \midrule
    \multicolumn{7}{l}{$\polP_1$}\\
    1\,604 & $1.35\times10^{-2}$ & -- & 52.3\% & $1.32\times10^{-2}$ & -- & 2.9\%\\
    5\,024 & $8.05\times10^{-3}$ & 0.91 & 58.9\% & $7.76\times10^{-3}$ & 0.93 & 2.2\%\\
    16\,023 & $4.82\times10^{-3}$ & 0.89 & 64.6\% & $4.53\times10^{-3}$ & 0.93 & 1.5\%\\
    49\,943 & $3.02\times10^{-3}$ & 0.82 & 66.6\% & $2.85\times10^{-3}$ & 0.82 & 1.0\%\\
    159\,720 & $1.84\times10^{-3}$ & 0.85 & 69.4\% & $1.73\times10^{-3}$ & 0.86 & 0.7\%\\
    499\,452 & $1.14\times10^{-3}$ & 0.85 & 71.0\% & $1.07\times10^{-3}$ & 0.84 & 0.5\%\\
    \midrule
    \multicolumn{7}{l}{$\polP_2$}\\
    1\,577 & $1.15\times10^{-2}$ & -- & 54.4\% & $1.05\times10^{-2}$ & -- & 5.7\%\\
    5\,161 & $7.20\times10^{-3}$ & 0.79 & 61.3\% & $6.64\times10^{-3}$ & 0.77 & 3.5\%\\
    16\,167 & $4.35\times10^{-3}$ & 0.88 & 64.6\% & $4.07\times10^{-3}$ & 0.86 & 2.6\%\\
    50\,333 & $2.87\times10^{-3}$ & 0.73 & 65.3\% & $2.71\times10^{-3}$ & 0.71 & 1.8\%\\
    159\,061 & $1.86\times10^{-3}$ & 0.75 & 67.3\% & $1.78\times10^{-3}$ & 0.74 & 1.5\%\\
    498\,969 & $1.21\times10^{-3}$ & 0.75 & 69.0\% & $1.19\times10^{-3}$ & 0.70 & 1.2\%\\
    \midrule
    \multicolumn{7}{l}{$\polP_3$}\\
    1\,678 & $1.10\times10^{-2}$ & -- & 43.4\% & $1.08\times10^{-2}$ & -- & 3.4\%\\
    5\,020 & $6.76\times10^{-3}$ & 0.89 & 51.0\% & $6.51\times10^{-3}$ & 0.92 & 2.5\%\\
    15\,409 & $4.30\times10^{-3}$ & 0.81 & 55.5\% & $4.20\times10^{-3}$ & 0.78 & 2.2\%\\
    49\,534 & $2.75\times10^{-3}$ & 0.76 & 56.5\% & $2.67\times10^{-3}$ & 0.78 & 2.1\%\\
    159\,115 & $1.74\times10^{-3}$ & 0.79 & 58.3\% & $1.71\times10^{-3}$ & 0.77 & 2.3\%\\
    497\,257 & $1.13\times10^{-3}$ & 0.76 & 59.2\% & $1.12\times10^{-3}$ & 0.73 & 1.8\%\\
    \midrule
    \multicolumn{7}{l}{$\polP_4$}\\
    1\,593 & $1.11\times10^{-2}$ & -- & 48.5\% & $1.06\times10^{-2}$ & -- & 5.9\%\\
    5\,045 & $6.84\times10^{-3}$ & 0.83 & 56.3\% & $6.59\times10^{-3}$ & 0.82 & 4.8\%\\
    16\,001 & $4.58\times10^{-3}$ & 0.69 & 59.5\% & $4.50\times10^{-3}$ & 0.66 & 3.8\%\\
    49\,561 & $3.03\times10^{-3}$ & 0.73 & 61.8\% & $3.06\times10^{-3}$ & 0.68 & 3.2\%\\
    161\,045 & $2.00\times10^{-3}$ & 0.70 & 63.0\% & $2.05\times10^{-3}$ & 0.68 & 2.7\%\\
    503\,301 & $1.29\times10^{-3}$ & 0.77 & 65.1\% & $1.34\times10^{-3}$ & 0.75 & 2.3\%\\
    \bottomrule
  \end{tabular}
\end{table}

\subsection{Mach 3 flow past a cylinder}
\begin{figure}[htbp]
  \centering
  \begin{minipage}[c]{0.54\textwidth}
    \begin{subfigure}{\textwidth}
      \centering
      \caption*{$\polP_1$, 999\,138 nodes}
      \includegraphics[width=\textwidth]{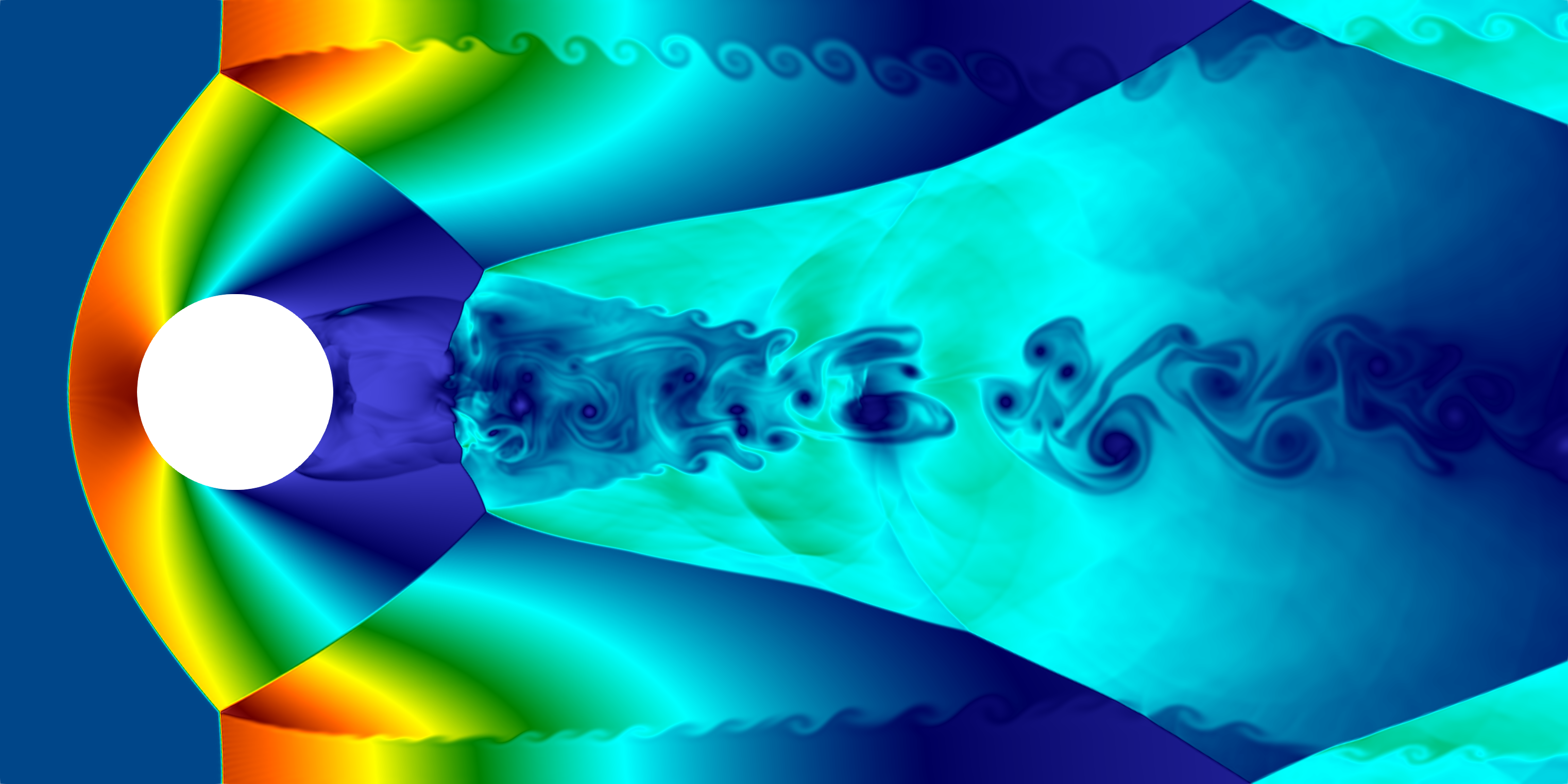}
    \end{subfigure}
    \begin{subfigure}{\textwidth}
      \centering
      \caption*{$\polP_2$, 999\,166 nodes}
      \includegraphics[width=\textwidth]{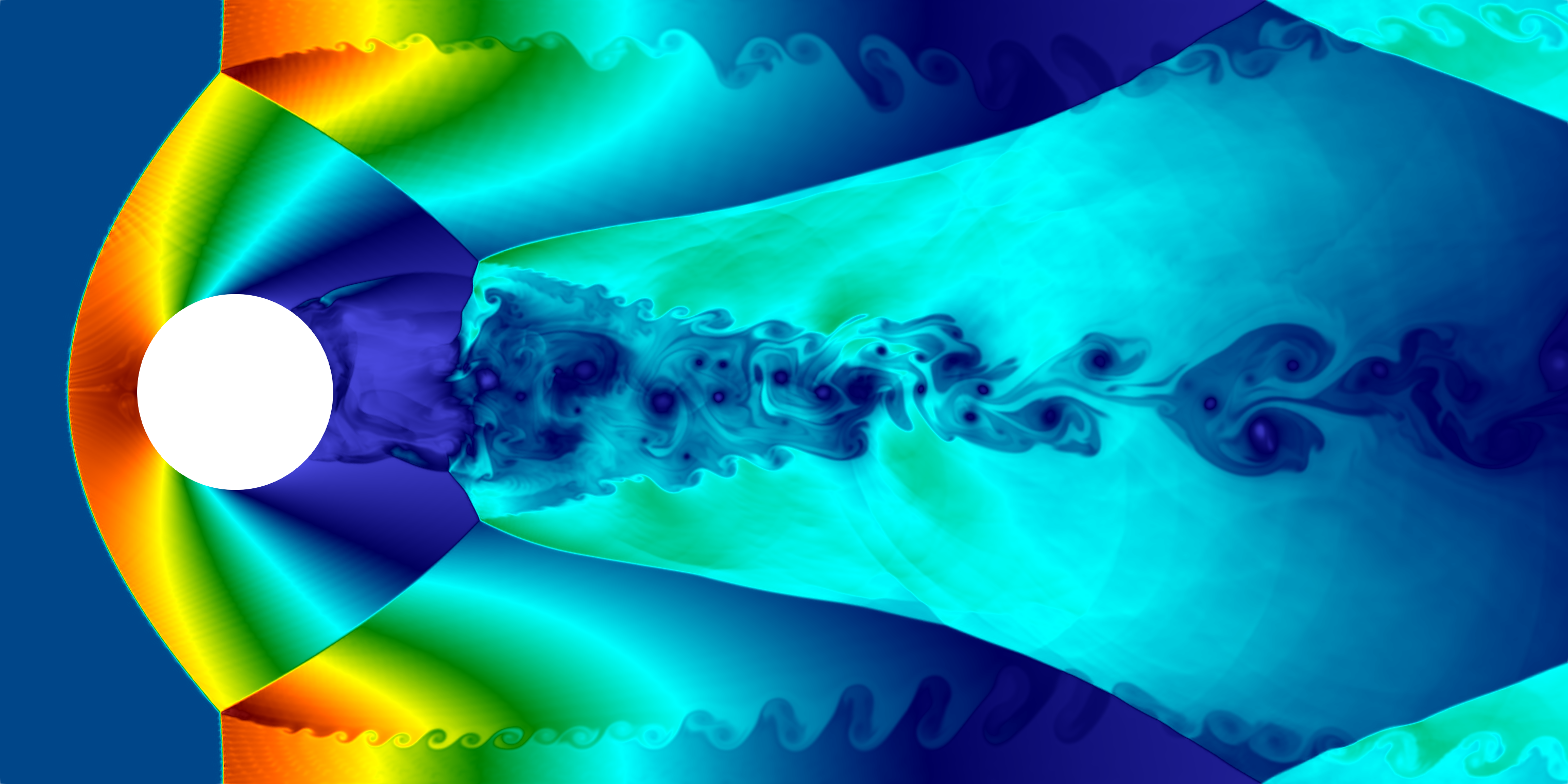}
    \end{subfigure}
    \begin{subfigure}{\textwidth}
      \centering
      \caption*{$\polP_3$, 996\,993 nodes}
      \includegraphics[width=\textwidth]{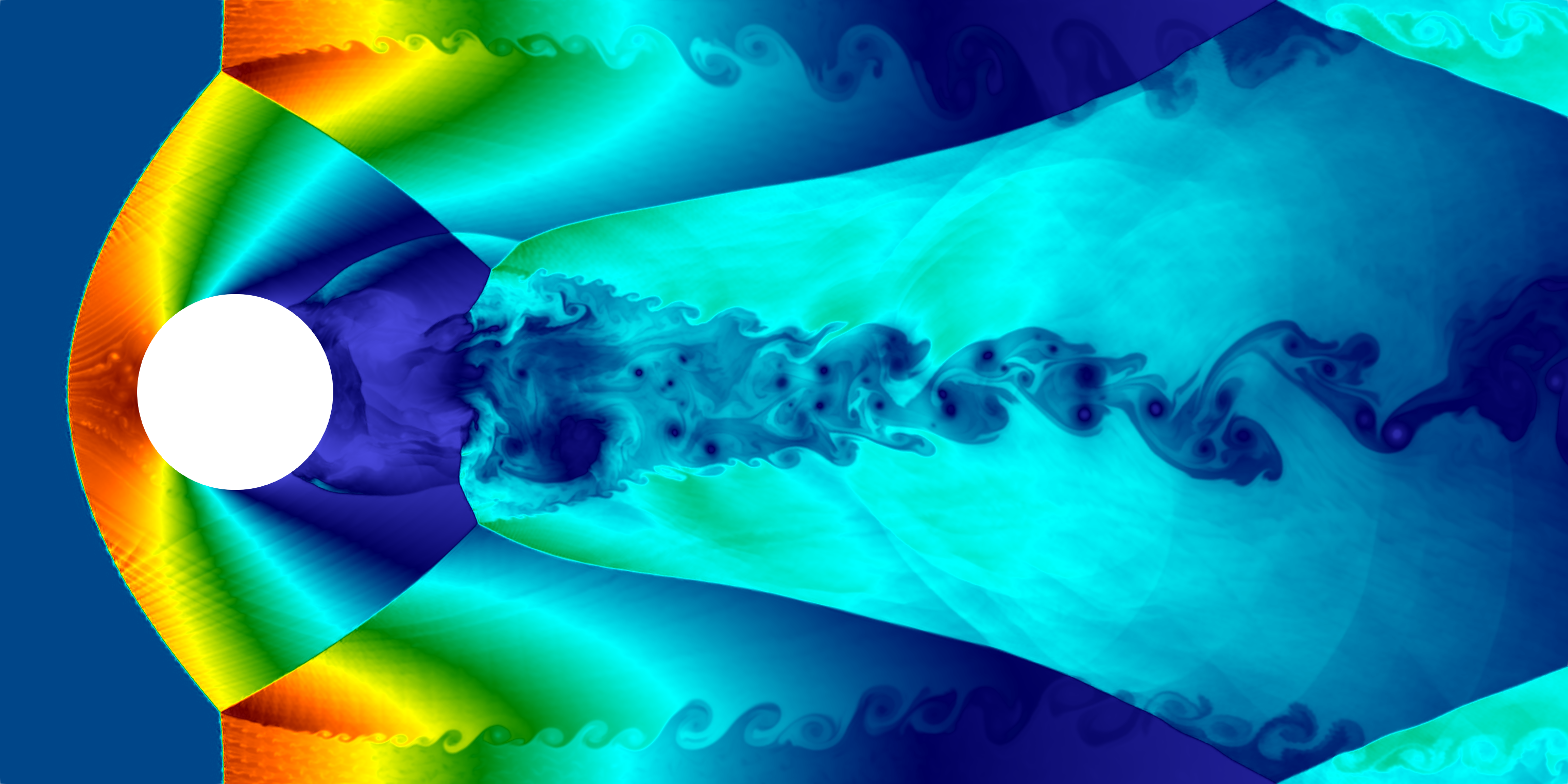}
    \end{subfigure}
    \begin{subfigure}{\textwidth}
      \centering
      \caption*{$\polP_4$, 995\,176 nodes}
      \includegraphics[width=\textwidth]{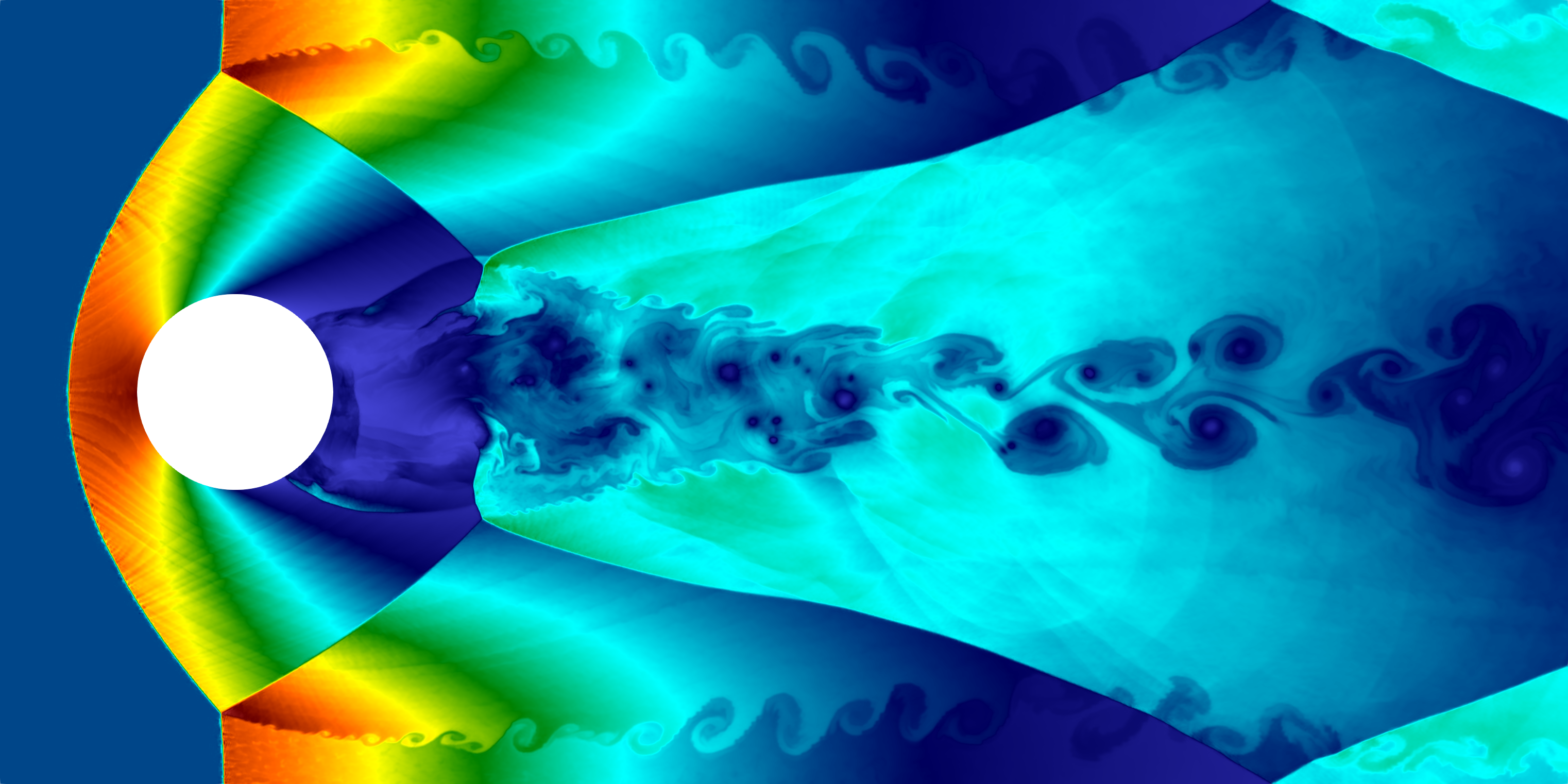}
    \end{subfigure}
  \end{minipage}%
  \hspace{0.02\textwidth}%
  \begin{minipage}[c]{0.1\textwidth}
    \includegraphics[width=\linewidth]{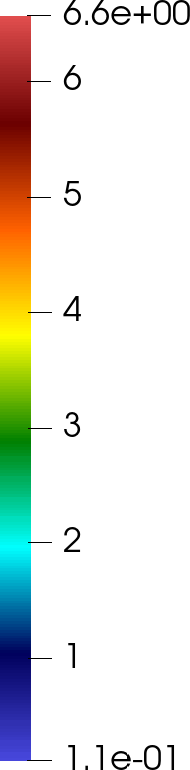}
  \end{minipage}
  \caption{
    Mach 3 flow past a cylinder at $\widehat t=4$ with element-local graph-Poisson fluxes, about $10^6$ nodes for each polynomial degree.}
  \label{fig:meshes}
\end{figure}

We apply the method to the Mach 3 flow past a cylinder of
\citep[Sec.~5.7]{Guermond_Nazarov_Popov_Tomas_2018}. The domain is
$
\Omega = \big( (0,4) \times (0,2) \big) \ \setminus \
\overline{B_{0.25}\big((0.6, 1) \big)},
$
where
$
B_{r}\big((x_0, y_0) \big) = \big\{ (x,y) \in \polR^2 \ | \ (x-x_0)^2 + (y - y_0)^2 < r^2 \big\}. 
$
The initial state and the inflow state at $x=0$ are $\rho = 1.4$, $\bv = (3,0)$, $p=1$. On the walls
and on the cylinder we impose slip-wall conditions, and no boundary condition is imposed
at the supersonic outflow at $x=4$. The meshes have about $10^6$ nodes for all polynomial degrees,
see Figure~\ref{fig:meshes}. We use SSP-RK(3,3). We also set
$C_{\textrm{R}} = 0.25$.

A strong bow shock forms in front of the cylinder and reflects from the walls, which produces Mach stems and triple points. Kelvin--Helmholtz instabilities develop along the slip lines that start at the triple points, and the small vortices in the wake of the cylinder are well resolved for all polynomial degrees.

\subsection{Mach 800 astrophysical jet}
\label{Sec:jet}
Finally, to test the limiter in an extreme case, we consider the Mach 800
astrophysical jet introduced by Balsara \citep{Balsara_2012}, in the setup of
\citep{Wu_Shu_2018} without the magnetic field, as in
\citep[Sec.~5.1.6]{Dzanic_Witherden_2022}, rotated by $90^\circ$ so that the
jet runs along the $x$-axis, with $\gamma=1.4$ on the half domain
$\Omega=(0,1.5)\times(0,0.5)$. The domain is initially filled with ambient gas
at rest, $(\rho,\bv,p)=(0.14,\boldsymbol 0,1)$. The jet state
$(\rho,\bv,p)=(1.4,(800,0),1)$ is prescribed on $x=0$, $y\le0.05$; its sound
speed is $1$, so the Mach number of the jet is $800$, and $253$ with respect to
the ambient sound speed. On the axis $y=0$ we impose the slip-wall (symmetry)
condition. At the nodes on $x=0$, $y>0.05$, we impose a far-field condition:
after every stage, the state is replaced by the solution at the boundary of the
one-dimensional Riemann problem in the normal direction between the low-order
state (inside) and the ambient gas at rest (outside), computed with the exact
Riemann solver of \citep[Ch.~4]{Toro_2009}. The nodes on $x=1.5$ and $y=0.5$
keep the low-order state. Both are treated as prescribed nodes, see
Section~\ref{Sec:boundary}.

We use meshes with about $1.09\cdot10^6$ nodes for every degree,
element-local fluxes, SSP-RK(3,3) with $\textrm{CFL}=0.3$ as in the other
examples, and $C_{\textrm{R}}=0.25$, until $\widehat t=0.002$.

\begin{figure}[htbp]
  \centering
  \begin{minipage}[c]{0.83\textwidth}
    \begin{subfigure}{\textwidth}
      \centering
      \caption*{$\polP_1$, 1\,085\,260 nodes}
      \includegraphics[width=\textwidth]{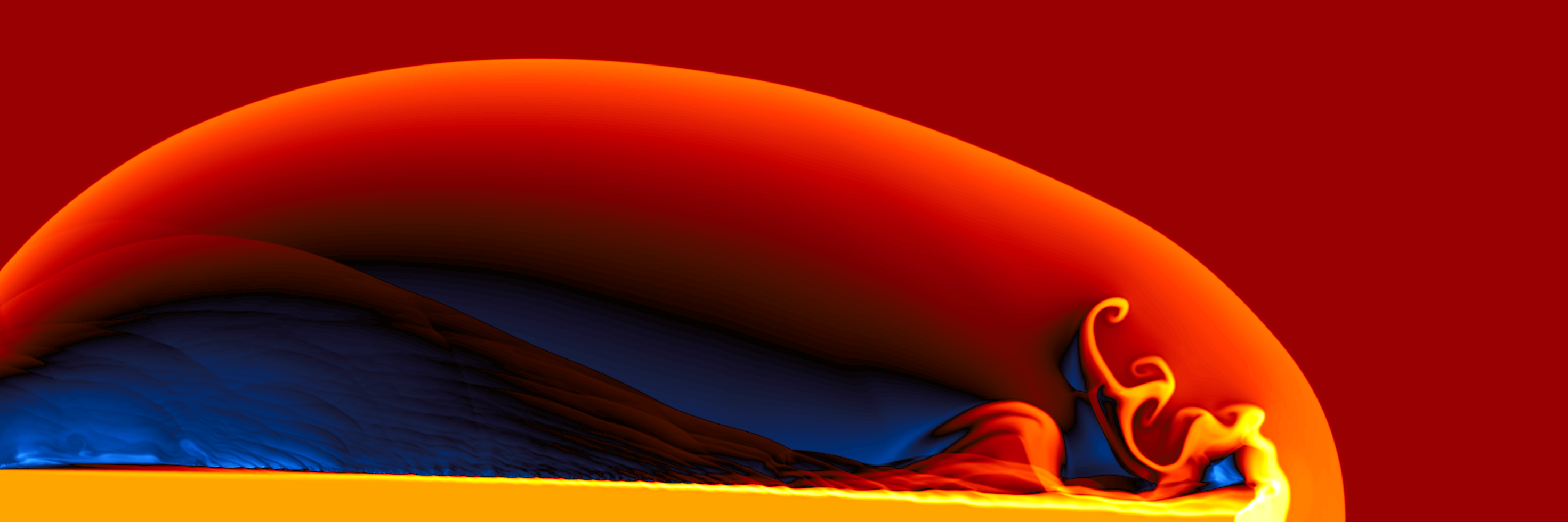}
    \end{subfigure}
    \begin{subfigure}{\textwidth}
      \centering
      \caption*{$\polP_2$, 1\,086\,269 nodes}
      \includegraphics[width=\textwidth]{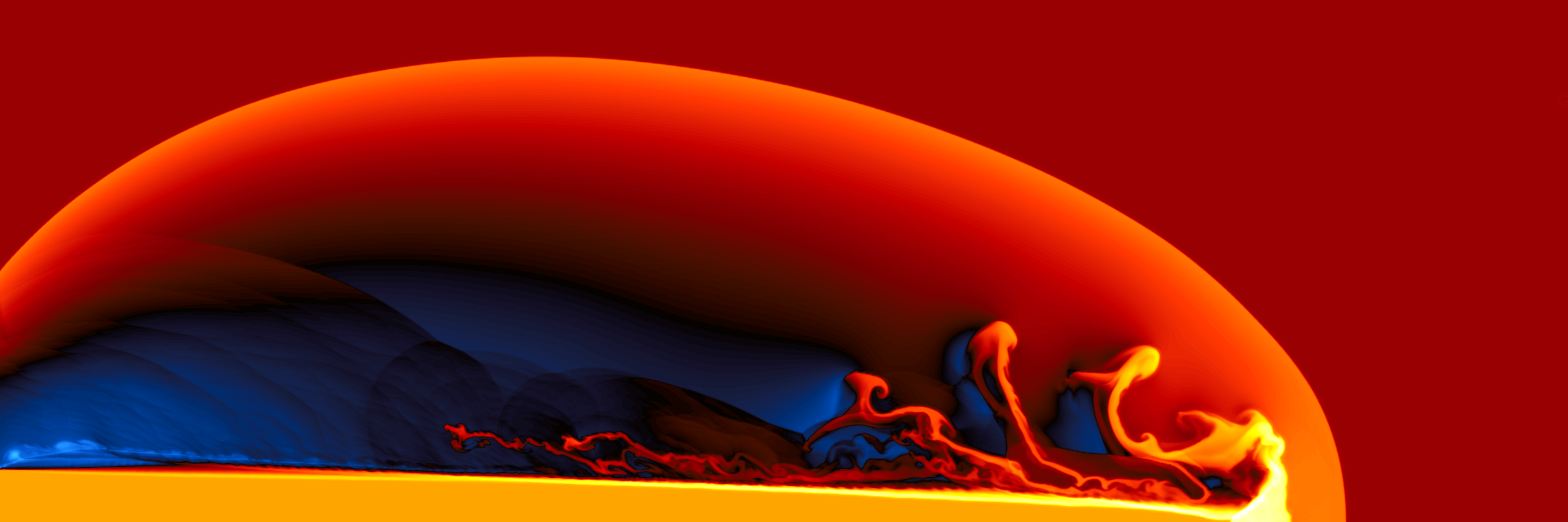}
    \end{subfigure}
    \begin{subfigure}{\textwidth}
      \centering
      \caption*{$\polP_3$, 1\,090\,600 nodes}
      \includegraphics[width=\textwidth]{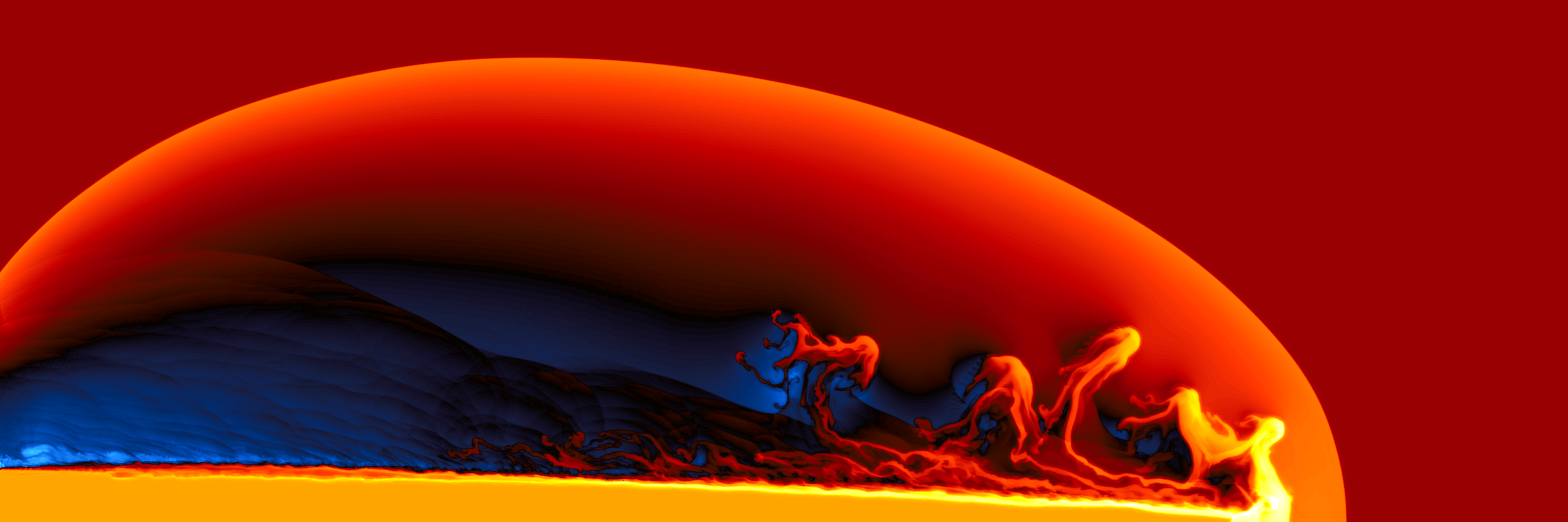}
    \end{subfigure}
    \begin{subfigure}{\textwidth}
      \centering
      \caption*{$\polP_4$, 1\,088\,541 nodes}
      \includegraphics[width=\textwidth]{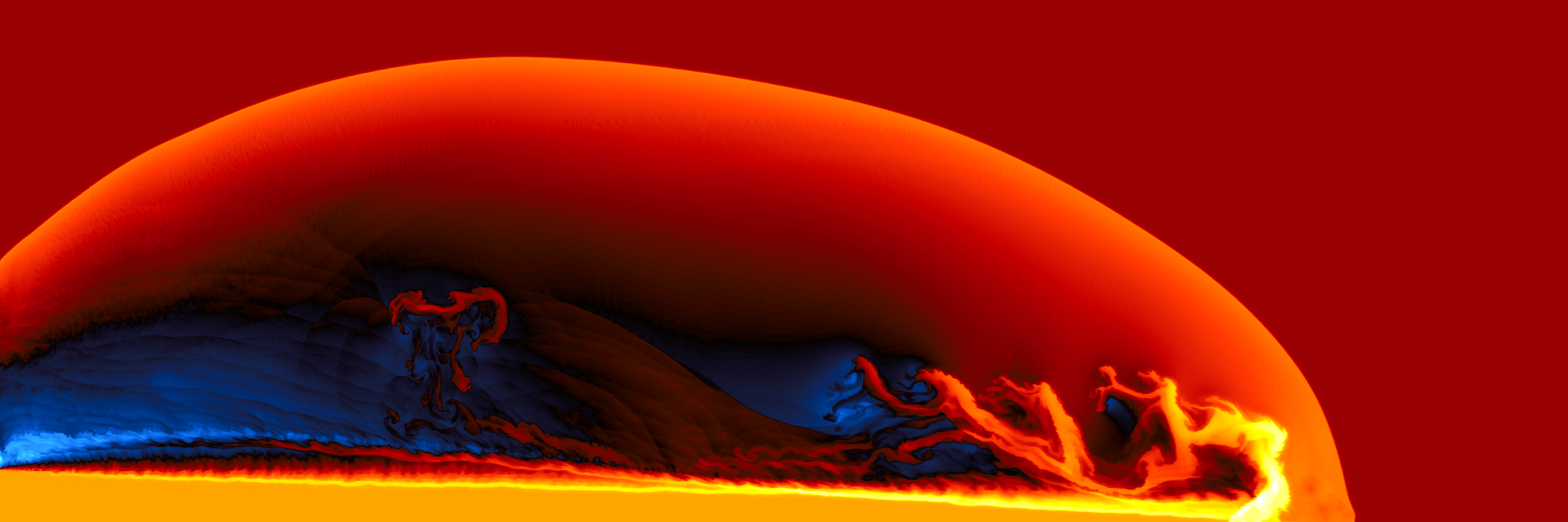}
    \end{subfigure}
  \end{minipage}%
  \hspace{0.02\textwidth}%
  \begin{minipage}[c]{0.1\textwidth}
    \includegraphics[width=\linewidth]{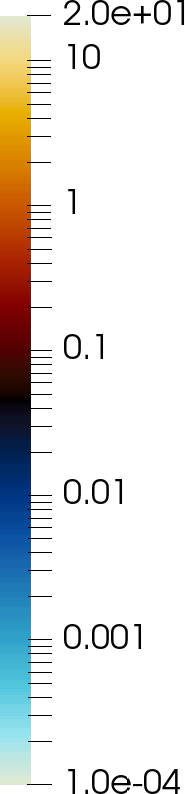}
  \end{minipage}
  \caption{Mach 800 jet at $\widehat t=0.002$: $\log_{10}\rho$ with element-local
    graph-Poisson fluxes.}
  \label{fig:jet}
\end{figure}

The density at the final time is shown in Figure~\ref{fig:jet}. The jet is
surrounded by a cocoon of shocked jet gas and by the shocked ambient gas behind
the bow shock; Kelvin--Helmholtz roll-ups develop along the contact between
them and at the head of the jet, and they get finer with the polynomial degree
at about the same number of nodes. The positions of the bow shock and of the
head of the jet agree with the $\polP_3$ results of
\citep{Dzanic_Witherden_2022} on their coarse mesh to within $0.015$ for all
degrees, see Table~\ref{tab:jet}. The runs take about $8\,200$ time steps for $\polP_1$, $\polP_2$ and
$\polP_4$, and about $12\,800$ for $\polP_3$, because of the factor
$r_{\max}=1.48$ in \eqref{eq:dt}.
The fraction of limited edges is $1.5\%$ for
$\polP_1$, $2.7\%$ for $\polP_2$, $2.8\%$ for $\polP_3$ and $4.7\%$ for
$\polP_4$, \ie $\polP_4$ limits almost twice as many edges as $\polP_2$ and
$\polP_3$. All runs reach the final time with positive
density and pressure, and the local bounds of
Section~\ref{Sec:convex_limiting} hold up to $1.5\cdot10^{-10}$, the tolerance
of the line search. The smallest density occurs at the far-field nodes next to
the nozzle, where the jet entrains the ambient gas.

\begin{table}[htbp]
  \caption{Mach 800 jet, $\widehat t=0.002$, about $1.09\cdot10^6$ nodes for every
    degree. Positions of the bow shock and of the head of the jet (where the
    density first deviates from the ambient value by more than $10\%$, and
    where $\rho>5$, respectively), compared with the $\polP_3$ results of
    \citep{Dzanic_Witherden_2022} on their coarse mesh (read off from their Fig.~6,
    $\pm0.005$); smallest and largest density at the final time;
    fraction of the edges whose flux is limited in the first pass.}
  \label{tab:jet}
  \centering
  \scriptsize
  \begin{tabular}{lccccc}
    \toprule
    & \citep{Dzanic_Witherden_2022} & $\polP_1$ & $\polP_2$ & $\polP_3$ & $\polP_4$\\
    \midrule
    \#nodes & & 1\,085\,260 & 1\,086\,269 & 1\,090\,600 & 1\,088\,541\\
    bow shock on $x=0$, $y$ & 0.245 & 0.242 & 0.242 & 0.242 & 0.254\\
    maximal height of the bow shock & 0.446 & 0.444 & 0.447 & 0.446 & 0.447\\
    bow shock on the axis, $x$ & 1.28 & 1.286 & 1.287 & 1.288 & 1.290\\
    end of the head ($\rho>5$), $x$ & 1.22 & 1.231 & 1.233 & 1.234 & 1.232\\
    $\min\rho$ & & $1.3\times10^{-3}$ & $5.7\times10^{-4}$ & $9.8\times10^{-5}$ & $1.0\times10^{-4}$\\
    $\max\rho$ & & 19.5 & 19.2 & 23.1 & 17.1\\
    time steps & & 8217 & 8165 & 12773 & 8228\\
    limited edges & & 1.5\% & 2.7\% & 2.8\% & 4.7\%\\
    \bottomrule
  \end{tabular}
\end{table}

\subsection{Efficiency comparison of global and local fluxes}
In this section, we compare the computational efficiency of two flux reconstructions:
global fluxes from Section~\ref{Sec:global}, and element-local fluxes from
Section~\ref{Sec:element_local}. We collected the computational time for
isentropic vortex of Section~\ref{Sec:vortex} in Table~\ref{tab:vortex_time} and those for
the Sod problem of Section~\ref{Sec:sod} in Table~\ref{tab:sod_time}.
All runs are serial, one core per run, on an Apple M2 Ultra; the runs
with the two flux reconstructions for the same degree were carried out
at the same time. 

Both reconstructions use the same time steps, and for the isentropic
vortex they give the same errors, to four digits, and limit almost no edges.
The difference in run time is therefore the cost of constructing the fluxes.
The element-local fluxes are faster for all degrees and meshes. On the finest
meshes, with about $5\cdot10^5$ nodes, the limiter is $1.4$--$3.5$ times
faster and the whole time stepping $1.2$--$2.1$ times faster for the
isentropic vortex; for the Sod problem the factors are $1.6$--$2.4$ and
$1.2$--$1.6$, where the global fluxes also limit many more edges, see
Table~\ref{tab:sod}. The gain grows with the number of nodes. The element
problems have the same size on every mesh, so their cost grows linearly with
number of nodes. The graph-Poisson problem couples all nodes; its sparse factorization is
computed once, but every solve costs more than $\calO(\#\textrm{nodes})$
because of the fill-in. In these serial runs, the graph Laplacian is factorized by the sparse LU factorization of PETSc.
The gain is largest for $k=1$, where the element-local fluxes are computed edge
by edge, without element problems.

\begin{table}[h!]
  \caption{Isentropic vortex, $\widehat t=2$: run times in seconds of the time stepping (serial, limiter in every stage) with the element-local (el) and the global (glob) graph-Poisson fluxes. The limiter time contains the flux construction: the element problems for el, the graph-Poisson solves for glob.}
  \label{tab:vortex_time}
  \centering
  \scriptsize
  \begin{tabular}{lrrrrrrrr}
    \toprule
    & & & \multicolumn{3}{c}{total} & \multicolumn{3}{c}{limiter}\\
    \cmidrule(lr){4-6}\cmidrule(lr){7-9}
    & $\#\textrm{nodes}$ & steps & el & glob & glob/el & el & glob & glob/el\\
    \midrule
    $\polP_1$ & 50\,464 & 184 & 46.8 & 64.6 & 1.38 & 19.9 & 41.0 & 2.07\\
     & 159\,544 & 320 & 262 & 383 & 1.46 & 109 & 252 & 2.31\\
     & 499\,990 & 580 & 1473 & 3028 & 2.06 & 591 & 2053 & 3.47\\
    \midrule
    $\polP_2$ & 50\,233 & 180 & 85.8 & 95.3 & 1.11 & 38.0 & 51.0 & 1.34\\
     & 159\,997 & 320 & 517 & 599 & 1.16 & 217 & 318 & 1.47\\
     & 501\,117 & 679 & 4140 & 5100 & 1.23 & 1715 & 2772 & 1.62\\
    \midrule
    $\polP_3$ & 50\,485 & 266 & 115 & 130 & 1.13 & 46.6 & 65.0 & 1.40\\
     & 159\,421 & 493 & 694 & 813 & 1.17 & 261 & 406 & 1.55\\
     & 495\,874 & 1086 & 5639 & 6874 & 1.22 & 2088 & 3708 & 1.78\\
    \midrule
    $\polP_4$ & 49\,665 & 179 & 86.8 & 94.4 & 1.09 & 39.2 & 49.1 & 1.25\\
     & 161\,481 & 348 & 567 & 644 & 1.13 & 245 & 337 & 1.37\\
     & 503\,969 & 607 & 3756 & 4768 & 1.27 & 1697 & 2449 & 1.44\\
    \bottomrule
  \end{tabular}
\end{table}

\begin{table}[h!]
  \caption{Sod problem, $\widehat t=0.2$: run times in seconds of the time stepping (serial) with the element-local (el) and the global (glob) graph-Poisson fluxes. The limiter time contains the flux construction: the element problems for el, the graph-Poisson solves for glob.}
  \label{tab:sod_time}
  \centering
  \scriptsize
  \begin{tabular}{lrrrrrrrr}
    \toprule
    & & & \multicolumn{3}{c}{total} & \multicolumn{3}{c}{limiter}\\
    \cmidrule(lr){4-6}\cmidrule(lr){7-9}
    & $\#\textrm{nodes}$ & steps & el & glob & glob/el & el & glob & glob/el\\
    \midrule
    $\polP_1$ & 49\,943 & 1102 & 220 & 307 & 1.39 & 108 & 216 & 2.01\\
     & 159\,720 & 1981 & 1280 & 2766 & 2.16 & 618 & 1986 & 3.21\\
     & 499\,452 & 3605 & 9149 & 14931 & 1.63 & 4390 & 10728 & 2.44\\
    \midrule
    $\polP_2$ & 50\,333 & 1108 & 366 & 432 & 1.18 & 174 & 262 & 1.50\\
     & 159\,061 & 2025 & 3237 & 3744 & 1.16 & 1481 & 2287 & 1.54\\
     & 498\,969 & 3686 & 16502 & 20892 & 1.27 & 7468 & 12726 & 1.70\\
    \midrule
    $\polP_3$ & 49\,534 & 1667 & 462 & 568 & 1.23 & 192 & 324 & 1.69\\
     & 159\,115 & 3070 & 3849 & 4827 & 1.25 & 1536 & 2584 & 1.68\\
     & 497\,257 & 5561 & 20505 & 28514 & 1.39 & 7718 & 16425 & 2.13\\
    \midrule
    $\polP_4$ & 49\,561 & 1111 & 378 & 430 & 1.14 & 187 & 254 & 1.36\\
     & 161\,045 & 2091 & 3403 & 3825 & 1.12 & 1782 & 2222 & 1.25\\
     & 503\,301 & 3734 & 17314 & 21007 & 1.21 & 7756 & 12526 & 1.62\\
    \bottomrule
  \end{tabular}
\end{table}

\section{Conclusions}
\label{Sec:conclusion}
In this paper, we presented high-order IDP continuous finite element
methods with Lagrange elements of degree one to four for the compressible
Euler equations. The low-order scheme uses $\polP_1$ elements on the fine mesh
whose vertices are the nodes of the high-order $\polP_k$ mesh. Its stencil
therefore consists of the nearest neighbours on the fine mesh, so that it adds
only the diffusion needed for the IDP property on this stencil. The high-order
scheme is the Galerkin method stabilized by a graph residual viscosity and a
continuous interior penalty. The difference between the two updates is written
as antisymmetric edge fluxes, either by a global graph-Poisson problem or
element by element; the element-local version needs no global linear solve.

Convex limiting of these fluxes gives a conservative method that keeps local
bounds on the density, the positivity of the internal energy and a local
minimum principle on the specific entropy. This also holds for $\polP_2$ and
$\polP_4$, whose lumped masses are zero or negative at some nodes; for these
elements, the high-order mass is restored within the bounds by an
element-local exchange, completed if needed by further rounds and a global
round. In our simulations, two rounds restore it exactly for smooth solutions
and up to a small remainder near strong shocks.
The numerical tests confirm near-optimal convergence rates for smooth
solutions and the robustness of the method for problems with strong shocks,
including a Mach 800 astrophysical jet.

\bibliographystyle{plainnat}
\bibliography{ref}
\end{document}